
%
%
\def\ifmagstep{\if00}

\input amstex
\loadbold
\documentstyle{amsppt}
\NoBlackBoxes

\ifmagstep
\pagewidth{32pc}
\pageheight{44pc}
\magnification=\magstep1
\else\fi

\def\ZZ{\Bbb{Z}}
\def\QQ{\Bbb{Q}}
\def\RR{\Bbb{R}}
\def\CC{\Bbb{C}}
\def\PP{\Bbb{P}}
\def\OO{\Cal{O}}
\def\Box{\square}
\def\Sym{\operatorname{Sym}}
\def\Gal{\operatorname{Gal}}
\def\Spec{\operatorname{Spec}}

\def\diag{\operatorname{diag}}
\def\Pic{\operatorname{Pic}}

\def\rank{\operatorname{rk}}
\def\id{\operatorname{id}}
\def\tr{\operatorname{tr}}
\def\ach{\widehat{\operatorname{ch}}}
\def\todd{\operatorname{td}}
\def\sch{\widetilde{\operatorname{ch}}}
\def\atodd{\widehat{\operatorname{td}}^A}
\def\chern{\operatorname{ch}}
\def\Ker{\operatorname{Ker}}
\def\Coker{\operatorname{Coker}}

\def\Ext{\operatorname{Ext}}
\def\Supp{\operatorname{Supp}}
\def\GL{\operatorname{GL}}
\def\rest#1#2{\left.{#1}\right\vert_{{#2}}}
\def\achern#1#2{\widehat{c}_{#1}(#2)}
\def\achsup#1#2{\overline{\operatorname{ch}}_{#1}(#2)}
\def\SStable{\overline{\boldkey M}}
\def\emotimes#1#2{\overset{#1}\to{\underset{#2}\to{\otimes}}}

\topmatter
\title
Inequality of Bogomolov-Gieseker's type \\ on arithmetic surfaces
\endtitle
\rightheadtext{}
\author Atsushi Moriwaki \endauthor
\leftheadtext{}
\affil
Department of Mathematics \\
University of California at Los Angeles \\
Los Angeles, CA 90024, USA \\
e-mail : moriwaki\@math.ucla.edu \\
April, 1993 (Revised Version)
\endaffil
\address
Department of Mathematics, Faculty of Science,
Kyoto University, Kyoto, 606-01, Japan
\endaddress
\curraddr
Department of Mathematics, University of California,
Los Angeles, 405 Hilgard Avenue, Los Angeles, California 90024, USA
\endcurraddr
\email moriwaki\@math.ucla.edu \endemail
\abstract
Let $K$ be an algebraic number field, $O_K$ the ring of integers of $K$,
and $f : X \to \Spec(O_K)$ an arithmetic surface.
Let $(E, h)$ be a rank $r$ Hermitian vector bundle on $X$
such that $E_{\overline{\QQ}}$ is
semistable on the geometric generic fiber
$X_{\overline{\QQ}}$ of $f$. In this paper, we will prove an arithmetic
analogy of Bogomolov-Gieseker's inequality:
$$
\achern{2}{E, h} - \frac{r-1}{2r} \achern{1}{E, h}^2 \geq 0.
$$
\vbox to 9cm {\vfil}
\endabstract
\toc\nofrills{Table of Contents}
\widestnumber\head{10}
\specialhead {} Introduction \endspecialhead
\head 1. Quick review of arithmetic intersection theory \endhead
\head 2. Asymptotic behavior of analytic torsions \endhead
\head 3. Hermitian modules over arithmetic curves \endhead
\head 4. Asymptotic behavior of $L^2$-degree of submodules of $H^0(L^m)$
\endhead
\head 5. Vanishing of a certain Bott-Chern secondary class \endhead
\head 6. Donaldson's Lagrangian and arithmetic second Chern class \endhead
\head 7. Second fundamental form \endhead
\head 8. Proof of the main theorem \endhead
\head 9. Arithmetic second Chern character of semistable vector bundles
\endhead
\head 10. Torsion vector bundles \endhead
\specialhead {} References \endspecialhead
\endtoc
\endtopmatter
\vfill\eject

\document

\head Introduction
\endhead

Let $M$ be an $n$-dimensional compact K\"{a}hler manifold with
a K\"{a}hler form $\Phi$. For a torsion free sheaf $F$ on $M$, we define
an averaged degree $\mu(F, \Phi)$ of $F$ with respect to $\Phi$ by
$$ \mu(F, \Phi) = \frac{{\displaystyle \int_M c_1(F) \wedge \Phi^{n-1}}}
{{\displaystyle \rank F}}.$$
Let $E$ be a torsion free sheaf on $M$. $E$ is said to be
{\it $\Phi$-stable} (resp. {\it $\Phi$-semisatble})
if, for every subsheaf $F$ of $E$ with $ 0 \subsetneq F \subsetneq E$,
an inequality $\mu(F, \Phi) < \mu(E, \Phi)$
(resp. $\mu(F, \Phi) \leq \mu(E, \Phi)$)
is satisfied.
If a torsion free sheaf $E$ of rank $r$ is $\Phi$-semistable, then we have
$$
        \int_M \left( c_2(E) - \frac{r- 1}{2 r} c_1(E)^2 \right)
        \wedge \Phi^{n-2} \geq 0,
$$
which is called Bogomolov-Gieseker's inequality
(cf. \cite{Bo} and \cite{Gi}).
The purpose of this paper is to establish an arithmetic analogy
of the above inequality, that is,

\proclaim{Main Theorem}
Let $K$ be an algebraic number field, $O_K$ the ring of integers of $K$,
and $f : X \to \Spec(O_K)$ an arithmetic surface.
Let $(E, h)$ be a Hermitian vector bundle on $X$.
If $E_{\overline{\QQ}}$ is semistable on the geometric generic fiber
$X_{\overline{\QQ}}$ of $f$, then we have an inequality
$$
\achern{2}{E, h} - \frac{r-1}{2r} \achern{1}{E, h}^2 \geq 0,
$$
where $r = \rank E$ and
$\achern{1}{E,h}$ and $\achern{2}{E,h}$ are arithmetic Chern classes
introduced by Gillet-Soul\'{e} \cite{GS90b}.
\endproclaim

The theory of stable vector bundles is closely related to
the Yang-Mills theory.
Let $E$ be a vector bundle on the compact K\"{a}hler manifold $M$.
$E$ is said to be
{\it $\Phi$-poly-stable} if there is a direct sum
$E = E_1 \oplus \cdots \oplus E_s$ such that $E_i$ is $\Phi$-stable for
all $1 \leq i \leq s$ and $\mu(E_1) = \cdots = \mu(E_s)$.
Moreover, let $h$ be a Hermitian metric of $E$.
$h$ is called {\it Einstein-Hermitian} (resp. {\it weakly Einstein-Hermitian})
if there is a constant (resp. $C^{\infty}$-function) $\varphi$
such that $\sqrt{-1} \Lambda K(E, h) = \varphi \id_{E}$, where
$K(E, h)$ is the curvature of $(E, h)$.
The fundamental theorem concerning Einstein-Hermitian metric is

\proclaim{Theorem A}
{\rm (cf. \cite{Do83}, \cite{Do85}, \cite{NS} and \cite{UY})}\quad
$E$ has an Einstein-Hermitian metric if and only if
$E$ is $\Phi$-poly-stable.
\endproclaim

The above theorem plays a crucial role for the proof of the main theorem.
Here we give a rough sketch of the proof of the main theorem,
which, I think, is a brief summary of this paper.

\subhead Step 1 \endsubhead
A poly-stable vector bundle is semistable, but a semistable vector bundle
is not necessarily poly-stable. This means that
a semistable bundle does not necessarily
have an Einstein-Hermitian metric. In general, a semistable vector bundle $F$
has a filtration:
$$
   0 = F_0 \subset F_1 \subset \cdots \subset F_{l-1} \subset F_l = F
$$
such that $F_i/F_{i-1}$ is stable for all $i$ and
$\mu(F_1/F_0) = \cdots = \mu(F_l/F_{l-1})$, which is
called {\it a Jordan-H\"{o}lder filtration of $F$}.
The first step of the proof of the main theorem is
a reduction to the case where $E$ is poly-stable on each infinite fiber.
Unfortunately, arithmetic Chern classes are sensitive for extensions, that is,
the Bott-Chern secondary characteristic appears in several formula for
an exact sequence of Hermitian vector bundles.
So we need an exact calculation of a certain secondary Bott-Chern
characteristic, which will be treated in \S7.
It follows a good comparison of second Chern classes (cf. Corollary~7.4).
Using this result, we can do our first step.

\subhead Step 2 \endsubhead
By the step 1, we may assume that
$E$ is poly-stable on each infinite fiber.
Thus by Theorem~A, $E$ has an Einstein-Hermitian metric $h_{EH}$
on each infinite fiber. Here we need a comparison of
$$
\achern{2}{E, h} - \frac{r - 1}{2r}\achern{1}{E, h}^2
\quad\hbox{and}\quad
\achern{2}{E, h_{EH}} - \frac{r - 1}{2r}\achern{1}{E, h_{EH}}^2.
$$
The following theorem guarantees that we may assume that
$h$ is Einstein-Hermitian.

\proclaim{Theorem C}{\rm (cf. Theorem~6.3)} \quad
Let $K$ be an algebraic number field and $O_K$ the ring of integers.
We denote by $K_{\infty}$ the set of
all embeddings of $K$ into $\CC$.
Let $f : X \longrightarrow \Spec(O_K)$ be an arithmetic variety with
$\dim X = d \geq 2$, and $(H, h_H)$ a Hermitian line bundle on $X$
such that, for each $\sigma \in K_{\infty}$,
$c_1(H_{\sigma}, h_{H_{\sigma}})$ gives a K\"{a}hler form $\Phi_{\sigma}$ on
an infinite fiber $X_{\sigma}$.
Let $E$ be a vector bundle of rank $r$ on $X$.
For a Hermitian metric $h$ of $E$, we set
$$ \Delta(E, h) = \left(
\achern{2}{E, h} - \frac{r - 1}{2r}\achern{1}{E, h}^2 \right)\cdot
\achern{1}{H, h_H}^{d-2}.$$
If $E$ is $\Phi_{\sigma}$-poly-stable on each infinite fiber $X_{\sigma}$,
then we have;
\roster
\item "(1)" the set
$\Delta = \{ \Delta(E, h) \mid \hbox{$h$ is a Hermitian metric of $E$} \}$
has the absolute minimal value.

\item "(2)" $\Delta(E, h_0)$ attaches the minimal value of $\Delta$
if and only if
$h_0$ is weakly Einstein-Hermitian on each infinite fiber.
\endroster
\endproclaim

\subhead Step 3 \endsubhead
Let $\pi : \PP(E) \longrightarrow X$ be the projective bundle of $E$ and
$\OO(1)$ the tautological line bundle on $\PP(E)$. We set
$L = \OO(r) \otimes \pi^*(\det E)^{-1}$.
The arithmetic analogy of
the Grothendieck relation~(1.9.1) and Lemma~5.1 implies that
$$
\hbox{
$(L^{r+1}) \leq 0$ $\Longrightarrow$
${\displaystyle \frac{r-1}{2r} \achern{1}{E,h}^2 \leq \achern{2}{E,h}}$.}
$$
Thus, it is sufficient to show that $(L^{r+1}) \leq 0$.

\subhead Step 4 \endsubhead
Let $N$ be an ample line bundle on $X$ with $\deg(N_K) > 2g(X_K) - 2$.
The arithmetic Riemann-Roch theorem
due to Gillet-Soul\'{e} \cite{GS92} shows us that
$$
   \chi_{L^2}(L^n \otimes \pi^*N) +
   \frac{1}{2} \tau((L^n \otimes \pi^*N)_{\infty})
   = \deg(\ach(L^n \otimes \pi^*N) \cdot \atodd(T_{\PP(E)/O_K})).
$$
Since we have
$$
   \lim_{n\to\infty}
   \frac{\deg(\ach(L^n \otimes \pi^*N) \cdot \atodd(T_{\PP(E)/O_K}))}{n^{r+1}}
   = \frac{1}{(r+1)!} (L^{r+1}),
$$
in order to get $(L^{r+1}) \leq 0$,
it is sufficient to see that
\roster
\item "(a)" $\tau((L^n \otimes \pi^*N)_{\infty}) \leq O(n^r \log n)$.

\item "(b)" $\chi_{L^2}(L^n \otimes \pi^*N) \leq O(n^r \log n)$.
\endroster

\subhead Step 5 \endsubhead
(a) is a consequences of the following theorem, which is
a generalization of \cite{BV} and Proposition 2.7.8 of \cite{Vo}
to the case of semipositive line bundles.

\proclaim{Theorem C}{\rm (cf. Theorem~2.1)}\quad
Let $X$ be a compact K\"{a}hler manifold of dimension $n$,
$(L, h)$ a Hermitian line bundle on $X$, and
$(E, h_E)$ a Hermitian vector bundle on $X$.
Let $\zeta_{q, d}$ be the zeta function of the Laplacian
$\Box_{q, d}$ on $A^{0, q}(L^{d} \otimes E)$.
Let $H_L$ be the Hermitian form corresponding to the curvature form
$K(L, h)$ of $L$ and $k$ a non-negative integer.
If $H_L(x)$ is positive semi-definite and $\rank H_L(x) \geq k$ for
all $x \in X$, then for $n - k < q \leq n$ there is a constant $C$ such that
$$
    | \zeta'_{q, d}(0) | \leq C d^n \log(d)
$$
for all $d \geq 0$.
\endproclaim

\subhead Step 6 \endsubhead
(b) can be derived from the following theorem.

\proclaim{Theorem D}{\rm (cf. Theorem~4.1)}\quad
Let $K$ be an algebraic number field, $O_K$ the ring of integers of $K$
and $f : X \longrightarrow \Spec(O_K)$ an arithmetic variety over $O_K$.
Let $(L, h)$ be a Hermitian line bundle on $X$,
$(E, h_E)$ a Hermitian vector bundle on $X$, and
$h_m$ a Hermitian metric of $H^0(X, L^m \otimes E)$ induced by
$h^m \otimes h_E$.
Then there is a constant $C$ such that
$$
    \frac{\deg_{L^2}(F, h_F)}{\rank F} \leq C m \log m
$$
for all $m > 1$ and all submodules $F$ of $H^0(X, L^m \otimes E)$,
where $h_F$ is the Hermitian metric induced by $h_m$.
\endproclaim

In \S9, we will consider an invariant of a semistable vector bundle arising
from the arithmetic second Chern character.
Let $f : X \to \Spec(O_K)$ be a regular arithmetic surface and
$E$ a vector bundle on $X$ such that
$\deg(E_K) = 0$ and $E_{\overline{\QQ}}$ is semistable.
We denote by $\operatorname{Herm}(E)$ the set of all Hermitian metrics of $E$.
Here we set
$$
   \achsup{2}{E} = \sup_{h \in \operatorname{Herm}(E)} \ach_2(E, h).
$$
As a consequence of Theorem~A, we can easily see that
$\achsup{2}{E}$ is a non-positive real number.
Let $\SStable_{X_K/K}(r, 0)$ be
the moduli scheme of semistable vector bundles
on $X_K$ with rank $r$ and degree $0$.
Our first question is

\definition{Question E}
Is there a relation between $\achsup{2}{E}$ and
a height of the class $[E_K]$ in $\SStable_{X_K/K}(r, 0)$?
\enddefinition

\noindent
For example, if $X_K$ is an elliptic curve,
we can give an answer (cf. Corollary~9.7).
For a curve of higher genus, it is still an open problem.
To solve the above question, we think that first of all
we must construct a certain canonical height of $\SStable_{X_K/K}(r, 0)$.
Moreover, we think that the theta divisor on $\SStable_{X_K/K}(r, 0)$
gives rise to a canonical height.

In \S10, we consider the equality condition of the inequality
of the main theorem. For a geometric case, this is closely related to
flat vector bundles. So we need an arithmetic analogy of flat vector bundles.
We give one candidate for arithmetic flatness, that is, torsion vector bundles.
(see Definition~10.6 for the definition of ``of torsion type''.)
Our question concerning torsion vector bundles is the following.

\definition{Question F}
Let $X$ and $E$ be as above.
If $\ach_2(E, h) = 0$, then is $E$ of torsion type?
\enddefinition

\noindent
A partial answer is given in Proposition~10.8.
Two questions are related to each other.
For example, if $r = 1$, $E$ satisfies some numerical conditions and
$h$ is Einstein-Hermitian, then
$\ach_2(E, h) = -[K : \QQ]\operatorname{Height}([E_K])$ (c.f. Lemma~9.2),
where $\operatorname{Height}$
is the N\'{e}ron-Tate height on $\Pic^0(X_K)$.
Thus, if $\ach_2(E, h) = 0$, $E_K$ is a torsion point of $\Pic^0(X_K)$.

\bigskip
Finally we would like to express hearty thanks to Prof. S. Zhang
for his valuable suggestions.
\vfill\eject

\head 1. Quick review of arithmetic intersection theory
\endhead

In this section, we would like to fix several notations of this paper,
which gives a quick review of arithmetic intersection theory.
Details will be found in \cite{GS90a}, \cite{GS90b}, \cite{Fa2} and
\cite{SABK}.

\subhead 1.1 Polynomial map associated with a formal power series $\phi$
\endsubhead
Let $\phi$ be a symmetric formal power series of $n$ variables over $\RR$,
$A$ an $\RR$-algebra, and $M_n(A)$ the algebra
of $(n \times n)$-matrices of $A$.
We denote by $\phi^{(k)}$ the homogeneous component of $\phi$ of degree $k$.
Then we can easily construct the unique polynomial map
$$
   \Phi^{(k)} : M_n(A) \longrightarrow A
$$
such that $\Phi^{(k)}$ is invariant under
conjugation by $\GL_n(A)$ and
$$
\Phi^{(k)} \left( \operatorname{diag}(\lambda_1, \ldots,  \lambda_n)
\right)
= \phi^{(k)}(\lambda_1, \ldots, \lambda_n).
$$
When $I$ is a nilpotent subalgebra, we may define
$$
    \Phi = \sum_{k} \Phi^{(k)} : M_n(I) \longrightarrow A.
$$
Standard examples of power series are the following:
\medskip\par
(Chern classes) \quad $c_i = s_i(T_1, \ldots, T_n)$, where $s_i$ is the
$i$-th elementary symmetric polynomial.
\par
(Total Chern class) \quad ${\displaystyle c = \sum_{i \geq 0} c_i}$.
\par
(Chern character) \quad ${\displaystyle \chern(T_1, \ldots, T_n) =
\sum_{i=1}^{n} \exp(T_i)}$.
\par
(Todd class) \quad ${\displaystyle \todd(T_1, \ldots, T_n)
= \prod_{i=1}^{n} \frac{T_i}{1 - \exp(-T_i)} }$.

\subhead 1.2 Chern characteristic class attached to $\phi$
\endsubhead
Let $M$ be a complex manifold and $A^{p,p}(M)$ the space of
complex $C^{\infty}$-forms of type $(p, p)$. We put
$$
   A(M) = \bigoplus_{p \geq 0} A^{p, p}(M) \quad\hbox{and}\quad
   \tilde{A}(M) = A(M)/(\operatorname{Im}\partial +
                        \operatorname{Im}\bar{\partial}).
$$
Here we set
${\displaystyle d^c = \frac{\partial - \bar{\partial}}{4 \pi \sqrt{-1}}}$.
Then we have
${\displaystyle d d^c = \frac{\sqrt{-1}}{2 \pi} \partial \bar{\partial}}$.
Let $(E, h)$ be a Hermitian vector bundle on $M$ of rank $n$,
$K(E, h)$ the curvature of $(E, h)$, and $\phi$ as above.
The {\it characteristic class} of $(E, h)$ attached to $\phi$ is defined by
$$
   \phi (E, h) = \Phi \left( \frac{\sqrt{-1}}{2 \pi} K(E, h) \right)
   \in A(M).
$$

\subhead 1.3 Bott-Chern secondary characteristic class
\endsubhead
Let $\Cal{E} : 0 \to (S, h') \to (E, h) \to (Q, h'') \to 0$
be an exact sequence of Hermitian vector bundles on $M$.
(Here $h'$ and $h''$ are not necessarily the metrics induced by $h$.)
We introduce the {\it Bott-Chern secondary characteristic class}
$\tilde{\phi}({\Cal{E}})$ of $\Cal{E}$ attached to $\phi$.
This class is characterized by the following three properties:
\roster
\item "(i)" $dd^c \tilde{\phi}({\Cal{E}}) =
\phi(E, h) - \phi((S, h')\oplus(Q, h''))$.

\item "(ii)" For every homomorphism $f : N \longrightarrow M$
of complex manifolds,
$$
\tilde{\phi}(f^*(\Cal{E})) = f^* \tilde{\phi}(\Cal{E}).
$$

\item "(iii)" If $(E, h) = (S, h') \oplus (Q, h'')$,
then $\tilde{\phi}(\Cal{E}) = 0$.
\endroster
(Note that the axiom (i) is different from \cite{BGS} or \cite{GS90b}.)

For a later purpose, we will show how to construct it.
We denote the homomorphism $S \longrightarrow E$ by $\alpha$. Let $\PP^1$ be
the projective line, $\OO(1)$ the tautological line bundle on $\PP^1$ and
$\sigma$ a section of $H^0(\PP^1, \OO(1))$ such that
$\sigma(\infty) = 0$. Let $p : M \times \PP^1 \longrightarrow M$ and
$q : M \times \PP^1 \longrightarrow \PP^1$ be the canonical projections.
We set
$$
  \tilde{E} = \Coker \left(
       p^*(\alpha) \oplus (\operatorname{id} \otimes \sigma) :
       p^*(S) \longrightarrow p^*(E) \oplus (p^*(S) \otimes q^*(\OO(1)))
                     \right).
$$
Note that $\rest{\tilde{E}}{M \times \{ 0 \}} \simeq E$ and
$\rest{\tilde{E}}{M \times \{ \infty \}} \simeq S \oplus Q$.
Here we give a Hermitian metric $\tilde{h}$ on $\tilde{E}$ such that
$\rest{(\tilde{E}, \tilde{h})}{M \times \{ 0 \}}$ and
$\rest{(\tilde{E}, \tilde{h})}{M \times \{ \infty \}}$ are
isometric to $(E, h)$ and $(S, h') \oplus (Q, h'')$ respectively. Then
$\tilde{\phi}(\Cal{E})$ is given as follows:
$$
   \tilde{\phi}(\Cal{E}) = \int_{\PP^1} \phi(\tilde{E}, \tilde{h}) \log |z|^2.
$$

\subhead 1.4 Green current
\endsubhead
Here we assume that the complex manifold $M$ is compact.
We denote by $\Cal{D}^{k, k}(M)$ the space of currents on $M$ of type $(k, k)$.
Let $Z$ be a cycle on $M$ of codimension $p$.
An element $g \in \Cal{D}^{p-1, p-1}(M)$ is called
a {\it green current} of $Z$ if
${\displaystyle  dd^c g + \delta_Z}$ is smooth.
The smooth form $dd^c g + \delta_Z$ is denoted by $\omega(g)$.
For example, let $D$ be a divisor
on $M$ and $h$ a Hermitian metric of the line bundle $\OO_M(D)$.
Let $s$ be a rational section of $\OO_M(D)$ such that $\operatorname{div}(s) =
D$.
Then it is well known that
$$
     - dd^c \log(h(s, s)) + \delta_D = c_1(\OO_M(D), h).
$$
Thus $- \log(h(s,s))$ is a green current of $D$.

\subhead 1.5 Analytic torsion
\endsubhead
Moreover we assume that $M$ is a K\"{a}hler manifold with a K\"{a}hler form
$\Omega$. Let $(E, h)$ be a Hermitian vector
bundle on $M$.
Let $\Box_q$ be the Laplacian on $A^{0,q}(E)$ and
$0 < \lambda_1 < \lambda_2 < \cdots $ eigenvalues of $\Box_q$. We set
$$
\zeta_q(s) = \sum_{i > 0} \lambda_i^{-s}.
$$
It is well known that $\zeta_q$ extends meromorphically to the complex plane
and is holomorphic at $s = 0$. We define
the {\it analytic torsion} $\tau(E, h)$ of $(E, h)$ by
$$
   \tau(E, h) = \sum_{q > 0} (-1)^q q \zeta_q'(0).
$$

\subhead 1.6 Arithmetic variety and arithmetic Chow group
\endsubhead
Let $K$ be an algebraic number field and $O_K$ the ring of integers of $K$.
We set $S = \Spec(O_K)$ and the set of $\CC$-valued points of $K$, that is,
$$
   \{ \sigma : K \hookrightarrow \CC \mid
   \hbox{$\sigma$ is an embedding of the field $K$ to $\CC$} \},
$$
is denoted by $S_{\infty}$ or $K_{\infty}$.
A projective and flat morphism of integral schemes $\pi : X \longrightarrow S$
is called an {\it arithmetic variety over $O_K$} if the generic fiber of $f$
is smooth.
For an objects $\operatorname{Ob}$ of $X$, we denote by
$\operatorname{Ob}_{\sigma}$
the pullback of $\operatorname{Ob}$ by an embedding $\sigma \in S_{\infty}$.
Let $Z^p(X)$ be a free abelian group generated by cycles of codimension $p$.
Let $\widehat{Z}^p(X)$ be a group of pairs $(Z, \sum_{\sigma} g_{\sigma})$
such that $Z \in Z^p(X)$ and $g_{\sigma}$ is a green current
of $Z_{\sigma}$ on $X_{\sigma}$.
Let $\widehat{\operatorname{CH}}^p(X)$ be the quotient group of
$\widehat{Z}^p(X)$
divided by the subgroup generated by the following elements:
\roster
\item "(a)" $(\operatorname{div}(f), \sum_{\sigma} - \log |f_{\sigma}|^2)$,
where $f$ is a rational function on some
subvariety $Y$ of codimension $p-1$ and $\log |f_{\sigma}|^2$
is the current defined by
$$
(\log |f_{\sigma}|^2)(\gamma) =
        \int_{Y_{\sigma}} (\log |f_{\sigma}|^2)\gamma.
$$

\item "(b)" $(0, \sum_{\sigma} \partial(\alpha_{\sigma}) +
                       \bar{\partial}(\beta_{\sigma}))$, where
$\alpha_{\sigma} \in \Cal{D}^{p-2, p-1}(X_{\sigma})$,
$\beta_{\sigma} \in \Cal{D}^{p-1, p-2}(X_{\sigma})$.
\endroster
We define three homomorphisms:
$$
\align
\omega &:  \widehat{\operatorname{CH}}^p(X) \longrightarrow
         \bigoplus_{\sigma} A^{p,p}(X_{\sigma}), \\
z &: \widehat{\operatorname{CH}}^p(X) \longrightarrow
     \operatorname{CH}^p(X), \\
a &: \sum_{\sigma} \tilde{A}^{p-1, p-1}(X_{\sigma}) \longrightarrow
    \widehat{\operatorname{CH}}^p(X)
\endalign
$$
by
$$
\align
   \omega(Z, \sum_{\sigma} g_{\sigma} ) & = \sum_{\sigma} \left(
     dd^c g_{\sigma} + \delta_{Z_{\sigma}} \right), \\
z(Z, \sum_{\sigma} g_{\sigma}) & = Z, \\
a(\sum_{\sigma} g_{\sigma}) & = (0, \sum_{\sigma} g_{\sigma})
\endalign
$$
respectively.
Moreover, $\deg : \widehat{\operatorname{CH}}^{\dim X}(X) \to \RR$
is given by
$$
  \deg(\sum_{P \in X} n_P P, \sum_{\sigma} g_{\sigma}) =
  \sum_{P \in X} n_P \log \#(O_X/m_P) +
  \frac{1}{2} \sum_{\sigma} \int_{X_{\sigma}} g_{\sigma},
$$
where $m_P$ is the maximal ideal of $P$.

\subhead 1.7 First Chern class of Hermitian vector bundle
\endsubhead
Let $E$ be a vector bundle on the arithmetic variety $X$.
We say $(E, h)$ is a {\it Hermitian vector bundle}
if, for all $\sigma \in S_{\infty}$, $E_{\sigma}$ is equipped with
a Hermitian metric $h_{\sigma}$.
Let $s_1, \ldots, s_n$ be rational sections of $E$ such that
$s_1, \ldots, s_n$ give a basis of $E$ at the generic point of $X$.
$\achern{1}{E, h}$ is defined by the class of
$$
    \left(
     \operatorname{div}(s_1 \wedge \cdots \wedge s_n),
     \sum_{\sigma \in S_{\infty}} - \log \det
        \pmatrix
          h_{\sigma}(s_1, s_1) &\cdots &h_{\sigma}(s_1, s_n) \\
          \vdots               &\ddots &\vdots               \\
          h_{\sigma}(s_n, s_1) &\cdots &h_{\sigma}(s_n, s_n)
        \endpmatrix
    \right).
$$

\subhead 1.8 Intersection of arithmetic cycles
\endsubhead
We set
$$
  \widehat{\operatorname{CH}}(X) = \bigoplus_{p \geq 0}
\widehat{\operatorname{CH}}^p(X).
$$
$\widehat{\operatorname{CH}}(X)_{\QQ}$ has a natural pairing
$$
\widehat{\operatorname{CH}}^p(X)_{\QQ} \otimes
\widehat{\operatorname{CH}}^q(X)_{\QQ}
\longrightarrow \widehat{\operatorname{CH}}^{p+q}(X)_{\QQ}
$$
defined by
$$
(Y, \sum_{\sigma} g_{Y_{\sigma}}) \cdot (Z, \sum_{\sigma} g_{Z_{\sigma}})
= (Y \cdot Z, \sum_{\sigma} g_{Y_{\sigma}} * g_{Z_{\sigma}}),
$$
where $g_{Y_{\sigma}} * g_{Z_{\sigma}} =
g_{Y_{\sigma}} \delta_{Z_{\sigma}} + \omega(g_{Y_{\sigma}}) g_{Z_{\sigma}}$.
In particular,
$$
 \achern{1}{L, h} \cdot (0, \sum_{\sigma} g_{\sigma})
= (0, \sum_{\sigma} c_1(L_{\sigma}, h_{\sigma})g_{\sigma}).
$$

\subhead 1.9 Arithmetic characteristic classes
\endsubhead
Let $\phi$ be a symmetric formal power series of $n$ variables
with real coefficients.
Let $(E, h)$ be a Hermitian vector bundle of rank $n$
on the arithmetic variety $X$.
We introduce the {\it characteristic class}
$$
      \hat{\phi}(E, h) \in \widehat{\operatorname{CH}}(X)_{\RR}
$$
attached to $\phi$. This is characterized by the following four axioms
(cf. Theorem 4.1 of \cite{GS90b}) :
\roster
\item "(I)" For every morphism $f : Y \longrightarrow X$
of arithmetic varieties,
$$
f^*(\hat{\phi}(E, h)) = \hat{\phi}(f^*(E, h)).
$$

\item "(II)" If $(E, h)$ is a direct summand of Hermitian line bundles, i.e.
$$
(E, h) = (L_1, h_1) \oplus \cdots \oplus (L_n, h_n),
$$
then we get
$$
   \hat{\phi}(E, h) = \phi(\achern{1}{L_1, h_1}, \ldots, \achern{1}{L_n, h_n}).
$$

\item "(III)" If we set
$$
   \phi(T_1 + T, \ldots, T_n + T) = \sum_{i \geq 0}
         \phi_i(T_1, \ldots, T_n) T^i,
$$
for a Hermitian line bundle $(L, h')$,
$$
   \hat{\phi}((E, h) \otimes (L, h')) = \sum_{i \geq 0}
       \hat{\phi}_i(E, h)\achern{1}{L, h'}^i.
$$

\item "(IV)" $\omega(\hat{\phi}(E, h)) =
\sum_{\sigma} \phi(E_{\sigma}, h_{\sigma})$.
\endroster
Let $\Cal{E} : 0 \to (S, h') \to (E, h) \to (Q, h'') \to 0$
be an exact sequence of Hermitian vector bundles on $X$.
An important property of $\hat{\phi}$ is
$$
   \hat{\phi}(E, h) - \hat{\phi}((S, h') \oplus (Q, h'')) =
   a(\sum_{\sigma} \tilde{\phi}(\Cal{E}_{\sigma})).
$$
For example, if we put
$$
  \hat{c}_t(E, h) = \sum_{i=0}^{\rank E}  (-1)^i \hat{c}_i(E, h) t^{n-i},
$$
we have
$$
   \hat{c}_t(E, h) - \hat{c}_t(S, h')\hat{c}_t(Q, h'') =
   a\left( \sum_{\sigma}\sum_{i=0}^{\rank E}(-1)^i
   \tilde{c}_i(\Cal{E}_{\sigma}) t^{n-i} \right).
$$
We apply this formula to the following special situation.
Let $f : \PP(E) \longrightarrow X$ be the projective bundle of $E$,
$Q = \OO_Y(1)$ the tautological line bundle of $\PP(E)$ and
$S$ the kernel of $f^*(E) \longrightarrow Q$.
We give the submetric $h'$ on $S$ and the quotient metric $h''$
on $Q$ induced by $f^* h$. Applying the above formula to the exact sequence
$$
   \Cal{E} : 0 \to (S, h') \to (f^* E, f^* h) \to (Q, h'') \to 0,
$$
we get
$$
   \hat{c}_t(f^*E, f^*h) - \hat{c}_t(S, h')\hat{c}_t(Q, h'') =
   a\left( \sum_{\sigma}\sum_{i=0}^{\rank E}(-1)^i
   \tilde{c}_i(\Cal{E}_{\sigma}) t^{n-i} \right),
$$
which implies
$$
\sum_{i=0}^{\rank E}
(-1)^i f^* \achern{i}{E, h} \achern{1}{Q, h''}^{n-i}
    = a \left(
     \sum_{\sigma}\sum_{i=0}^{\rank E}(-1)^i
     \tilde{c}_i(\Cal{E}_{\sigma}) c_1(Q_{\sigma}, h''_{\sigma})^{n-i}.
        \right),
\tag 1.9.1
$$
by evaluating $t = \achern{1}{Q, h''}$.
This is an arithmetical analogy of the Grothendieck relation.
Conversely, the relation (1.9.1) defines
$\achern{i}{E, h}$. To see this, we will prove the following fact.
\proclaim{Claim 1.9.2}
\roster
\item "(1)" A homomorphism
$\psi : \widehat{\operatorname{CH}}(X)^{\oplus n} \longrightarrow
        \widehat{\operatorname{CH}}(Y)$ by
$$
\psi(x_0, \cdots, x_{n-1}) = \sum_{i=0}^{n-1} f^*(x_i) \achern{1}{Q, h''}^i
$$
is injective, where $n = \rank E$ and $Y = \PP(E)$.

\item "(2)" The image of $\psi$ is
$$
   \{ x \in \widehat{\operatorname{CH}}(Y) \mid
   \omega(x) \in \sum_{\sigma} \sum_{i=0}^{n-1}
   f_{\sigma}^*(A(X_{\sigma}))c_1(Q_{\sigma}, h''_{\sigma})^i \}.
$$
\endroster
\endproclaim

\demo{Proof}
(1) can be easily checked by the formula:
$$
f_*(f^*(\gamma) \achern{1}{Q, h''}^i) =
\cases
0, & 0 \leq i < n-1, \\
\gamma, & i = n-1.
\endcases
$$
Next we consider (2).
Let $x$ be an element of $\widehat{\operatorname{CH}}(Y)$ such that
$$
   \omega(x) \in \sum_{\sigma} \sum_{i=0}^{n-1}
   f_{\sigma}^*(A(X_{\sigma}))c_1(Q_{\sigma}, h''_{\sigma})^i.
$$
Since $\psi$ is bijective on finite part, we may assume that $z(x) = 0$.
Thus we can set $x = (0, \sum_{\sigma} g_{\sigma})$ with
$g_{\sigma} \in A(Y_{\sigma})$.
Here we put
$$
\omega(x) = \sum_{\sigma} \sum_{i=0}^{n-1}
f^*(u_{\sigma, i}) c_1(Q_{\sigma}, h''_{\sigma})^i.
$$
Then
$$
  dd^c (g_{\sigma}) =
  \sum_{i=0}^{n-1} f^*(u_{\sigma, i}) c_1(Q_{\sigma}, h''_{\sigma})^i.
$$
Using integration along $f_{\sigma}$, we can find $v_{\sigma, i}$ with
$dd^c v_{\sigma, i} = u_{\sigma, i}$.
(For example
${\displaystyle v_{\sigma, n-1} = \int_{f_{\sigma}} g_{\sigma}}$.)
Then
$$
\omega\left( x - \sum_{i=0}^{n-1}
f^*(0, \sum_{\sigma} v_{\sigma, i}) \achern{1}{Q, h''}^i \right) = 0.
$$
Thus we may assume that $\omega(x) = 0$,
which implies $g_{\sigma}$ is harmonic up to
$\operatorname{Im}\partial + \operatorname{Im}\bar{\partial}$.
Hence by the structure of cohomology ring $H^*(Y_{\sigma}, \RR)$,
$g_{\sigma}$ can be written as the form
$$
   \sum_{i=0}^{n-1} g_{\sigma, i} c_1(Q_{\sigma}, h''_{\sigma})^i
$$
up to $\operatorname{Im}\partial + \operatorname{Im}\bar{\partial}$.
Therefore
${\displaystyle x =
\sum_{i=0}^n f^*(0, \sum_{\sigma} g_{\sigma, i}) \achern{1}{Q, h''}^i}$.
\enddemo

\bigskip
Let us go back to the construction of $\achern{i}{E, h}$.
By this claim, there is a unique sequence
$\{ c_1, \ldots, c_n \}$ of $\widehat{\operatorname{CH}}(X)$ such that
$$
\sum_{i=1}^{n} (-1)^i f^* (c_i) \achern{1}{Q, h''}^{n-i} =
a\left( \sum_{\sigma}\sum_{i=0}^{n}(-1)^i
\tilde{c}_i(\Cal{E}_{\sigma}) c_1(Q_{\sigma}, h''_{\sigma})^{n-i} \right)
- \achern{1}{E, h}^n.
$$
Hence we may set $\achern{i}{E, h} = c_i$.

\subhead 1.10 $L^2$-degree of Hermitian module
\endsubhead
Let $V$ be $O_K$-module of finite rank.
A pair $(V, h)$ is called {\it a Hermitian module}
if we give a Hermitian metric $h_{\sigma}$ on $V_{\sigma}$ for each
$\sigma \in K_{\infty}$.
For example, $O_K$ has the canonical metric $\operatorname{can}_K$
induced by each embedding $\sigma : K \hookrightarrow \CC$.
We define $L^2$-degree
$\deg_{L^2}(V, h)$ of $(V, h)$ by
$$
\multline
  \deg_{L^2}(V, h) = \log\#\left( \frac{V}{O_K x_1 + \cdots + O_K x_t} \right)
\\
   - \frac{1}{2} \sum_{\sigma \in S_{\infty}} \log \det
   \pmatrix
   h_{\sigma}(x_1, x_1) & \cdots & h_{\sigma}(x_1, x_t) \\
   \vdots               & \ddots & \vdots               \\
   h_{\sigma}(x_t, x_1) & \cdots & h_{\sigma}(x_t, x_t) \\
   \endpmatrix ,
\endmultline
$$
where $x_1, \ldots, x_t \in M$ and $\{ x_1, \ldots, x_t \}$ is a basis
of $M \otimes K$. Using the Hasse product formula,
it is easily checked that $\deg_{L^2}(M, h)$ does not depend on the choice of
$\{ x_1, \ldots, x_t \}$.

\subhead 1.11 Arithmetic Riemann-Roch Theorem
\endsubhead
Finally we explain the arithmetic Riemann-Roch theorem.
Let $\pi : X \longrightarrow S$ an arithmetic variety.
We give a K\"{a}hler form $\Omega_{\sigma}$ for each
infinite fiber $X_{\sigma}$. Let $(E, h)$ be a Hermitian vector bundle on $X$.
Using this K\"{a}hler form,
we can give a metric of $H^q(X_{\sigma}, E_{\sigma})$
by identifying $H^q(M_{\sigma}, E_{\sigma})$ with
the space of harmonic forms ${\boldkey H}^{0,q}(E_{\sigma})$.
Hence we can define $\deg_{L^2}(H^q(X, E))$. Thus we have
$L^2$-Euler character
$$
 \chi_{L^2}(E, h) = \sum_{q \geq 0} (-1)^q \deg_{L^2}(H^q(X, E)).
$$
To state arithmetic Riemann-Roch, we need to mention about the arithmetic
Todd character.
Let $\zeta(s) = \sum_{n > 0} n^{-s}$ be the standard zeta function and
we define a formal power series $R(T_1, \ldots, T_n)$ as
$$
 \sum\Sb m \geq 1 \\ m : \operatorname{odd} \endSb
 \frac{1}{m!} \left(
    2 \zeta'(-m) + \zeta(-m)( 1 + 1/2 + 1/3 + \cdots + 1/m )
    \right) (T_1^m + \cdots + T_n^m).
$$
The arithmetic Todd character is defined by
$$
\atodd(T_{X/S}, \Omega) = \widehat{\operatorname{td}}(T_{Y/S}, \Omega)
                   (1 - a(\sum_{\sigma} R(T_{X_{\sigma}}, \Omega_{\sigma}))).
$$
The arithmetic Riemann-Roch due to Gillet-Soul\'{e} \cite{GS92}
is as follows:
$$
   \chi_{L^2}(E, h) + \frac{1}{2}
   \sum_{\sigma \in S_{\infty}}\tau(E_{\sigma}, h_{\sigma})
   = \deg(\ach(E, h) \cdot \atodd(T_{Y/S}, \Omega)).
\tag 1.11.1
$$
\vfill\eject

\head 2. Asymptotic behavior of analytic torsions
\endhead

In this section, we will consider asymptotic behavior of analytic torsions
of powers of a semi-positive line bundle.
The main theorem of this section is the following theorem.

\proclaim{Theorem 2.1}
Let $X$ be a compact K\"{a}hler manifold of dimension $n$,
$(L, h)$ a Hermitian line bundle on $X$, and
$(E, h_E)$ a Hermitian vector bundle on $X$.
Let $\zeta_{q, d}$ be the zeta function of the Laplacian
$\Box_{q, d}$ on $A^{0, q}(L^{d} \otimes E)$.
Let $H_L$ be the Hermitian form corresponding to the curvature form
$K(L, h)$ of $L$ and $k$ a non-negative integer.
If $H_L(x)$ is positive semi-definite and $\rank H_L(x) \geq k$ for
all $x \in X$, then for $n - k < q \leq n$ there is a constant $C$ such that
$$
    | \zeta'_{q, d}(0) | \leq C d^n \log(d)
$$
for all $d \geq 0$.
\endproclaim

First we prepare the following lemma.

\proclaim{Lemma 2.2}
With notation as in Theorem~{\rm 2.1},
if $H_L(x)$ is positive semi-definite and $\rank H_L(x) \geq k$ for
all $x \in X$, then, for $n - k < q \leq n$, there is a positive constant
$c_{L, q}$ such that $c_{L, q}$ depends only on $(L, h)$ and $q$, and that
$$
 - \sqrt{-1} \langle (\Lambda e(K_L) - e(K_L) \Lambda) \phi, \phi \rangle(x)
 \geq c_{L, q} \langle  \phi, \phi \rangle(x)
$$
for all $\phi \in A^{n, q}(E)$ and $x \in X$.
In particular, integrating the above inequality,
we have
$$
 - \sqrt{-1} \left( (\Lambda e(K_L) - e(K_L) \Lambda) \phi, \phi \right)
 \geq c_{L, q} ( \phi, \phi ).
$$
\endproclaim

\demo{Proof}
First we claim:

\proclaim{Claim 2.2.1}
There is a positive constant $c$ such that, for all $x \in X$,
we can take an orthogonal basis $\{ w_1, \ldots, w_n \}$ of $T_{X, x}$
with
$$
\cases
H_L(x)(w_i, w_j) = h_j \delta_{ij} & \\
h_i \geq c \quad (1 \leq i \leq k).
\endcases
$$
\endproclaim

For a point $p \in X$, there is an open neighborhood $U_p$ of $p$ and
a $C^{\infty}$-subvector bundle $F_p$ of $T_{U_p}$ such that
$\rank F_p = k$ and $\rest{H_L}{F_p}$ is positive definite.
Shrinking $U_p$ if necessarily, we may assume that
there is a positive constant $c_p$ such that $H_L(y)(v, v) \geq c_p$
for all $y \in U_p$ and all $v \in (F_p)(y)$ with $||v|| = 1$.
Since $X$ is compact, we have $p_1, \ldots, p_s \in X$ with
$\bigcup_{t=1}^{s} U_{p_t} = X$. We set
$c = \min \{c_{p_1}, \ldots, c_{p_s} \}$.
Let $x$ be an arbitrary point of $X$. If $x \in U_{p_t}$, then
we have an orthogonal basis $\{w_1, \ldots, w_n \}$ of $T_{X, x}$
such that $w_1, \ldots, w_k \in (F_{p_t})(x)$ and
$H_L(x)(w_i, w_j) = h_i\delta_{ij}$. Thus $h_i \geq c_{p_t} \geq c$
for $1 \leq i \leq k$.
Therefore we get our claim.

\medskip
Let $\{ w_1, \ldots, w_n \}$ be an orthogonal basis of $T_{X, x}$ as
in Claim~2.2.1 and $\{ \theta_1, \ldots, \theta_n \}$
the dual basis of $\{ w_1, \ldots, w_n \}$. We set
$$
 \phi = \sum_{i_1 < \cdots < i_q} \phi_{i_1, \ldots, i_q}
 \theta_1 \wedge\cdots\wedge \theta_n \wedge
 \bar{\theta}_{i_1}\wedge\cdots\wedge\bar{\theta}_{i_q}.
$$
Then by the formula due to Gigante (cf. p.70, (3.6) in \cite{Ko}),
we have
$$
\align
 - \sqrt{-1} \langle (\Lambda e(K_L) - e(K_L) \Lambda) \phi, \phi \rangle(x)
 & = \sum_{i_1 < \cdots < i_q}
       (h_{i_1} + \cdots + h_{i_q}) ||\phi_{i_1, \ldots, i_q}||_E^2(x) \\
 & \geq \sum_{i_1 < \cdots < i_q}
       (q - n + k)c ||\phi_{i_1, \ldots, i_q}||_E^2(x) \\
 & = (q - n + k) c \langle \phi , \phi \rangle(x).
\endalign
$$
Thus we obtain our lemma by setting $c_{L, q} = (q - n + k)c$.
\qed
\enddemo

\subhead 2.3 \endsubhead
Let us start the proof of Theorem~2.1.
The idea of this proof is found in \cite{BV} and \cite{Vo}.
Let $\lambda_{q, d}$ be the minimal eigenvalue of the Laplacian
$\Box_{q, d}$. By Lemma 2.7.7 of \cite{Vo},
it is sufficient to see that there are
constants $c > 0$ and $c' \geq 0$ such that
$$
     \lambda_{q, d} \geq cd - c'.
$$
Let $\Box'_{q, d}$ be the Laplacian on
$A^{n, q}(L^d \otimes E \otimes \omega_X^{-1})$.
As mentioned in the proof of Theorem 1 of \cite{BV},
$\Box_{q, d}$ and $\Box'_{q, d}$ have
the same spectrum. We set $E' = E \otimes \omega_X^{-1}$.
Let $K_{E'}$ and $K_d$ be curvature forms of $E'$ and
$L^{d} \otimes E'$ respectively. Then we have
$$
  K_d = d K_L + \operatorname{id} \otimes K_{E'}.
$$
For $\phi \in A^{n, q}(L^d \otimes E')$, by an easy calculation, we get
$$
\align
  (\Box'_{q, d} \phi, \phi) & =
  ||D' \phi||^2 + ||\delta' \phi||^2  -
  \sqrt{-1} \left( (\Lambda e(K_d) - e(K_d)\Lambda)\phi, \phi \right) \\
  & \geq
  - d \sqrt{-1}\left( (\Lambda e(K_L) - e(K_L)\Lambda)\phi, \phi \right) \\
  & \phantom{\geq - d \sqrt{-1}}
  -\sqrt{-1}\left( (\Lambda (\operatorname{id} \otimes e(K_{E'})) -
                 (\operatorname{id} \otimes e(K_{E'}))\Lambda)\phi, \phi
\right),
\endalign
$$
where $D'$ is $(1, 0)$-part of the Hermitian connection of $L^d \otimes E'$
and $\delta'$ is the adjoint operator of $D'$.
By virtue of Lemma~2.2,
there is a positive constant $c$ such that
$$
- \sqrt{-1}\left( (\Lambda e(K_L) - e(K_L)\Lambda)\phi, \phi \right)
 \geq c(\phi, \phi)
$$
for all $\phi \in A^{q,n}(L^d \otimes E')$.
On the other hand, we can easily find $c' \geq 0$ such that
$$
  - \sqrt{-1}\left( (\Lambda (\operatorname{id} \otimes e(R_{E'})) -
                 (\operatorname{id} \otimes e(R_{E'}))\Lambda)\phi, \phi
\right)
  \geq  -c' (\phi, \phi)
$$
for all $\phi \in A^{q,n}(L^d \otimes E')$.
Therefore we have an inequality
$$
  (\Box'_{d, q} \phi, \phi) \geq (cd - c')(\phi, \phi).
$$
We choose $\phi \in A^{q,n}(L^d \otimes E')$ such that
$\Box'_{q, d} \phi = \lambda_{q, d} \phi$ and $(\phi, \phi) = 1$.
Thus using the above inequality, we obtain
$$
    \lambda_{q, n} \geq cd - c'.
\eqno{\qed}
$$

\proclaim{Corollary 2.4}
Let $X$ be a compact K\"{a}hler manifold of dimension $n$,
$(L, h)$ a Hermitian line bundle on $X$, and
$(E, h_E)$ a Hermitian vector bundle on $X$.
For the Hermitian form $H_L$ corresponding to the curvature form of $(L, h)$,
we assume that
$H_L(x)$ is positive semi-definite and $\rank H_L(x) \geq n-1$ for
all $x \in X$.
Then there is a positive constant $C$ such that
$$
    \tau(L^d \otimes E) \leq C d^n \log(d)
$$
for all $d \geq 0$.
\endproclaim

\demo{Proof}
Let $\zeta_{q, d}$ be the zeta function of the Laplacian on
$A^{0,q}(L^d \otimes E)$. Then
$$
  \tau (L^d \otimes E) =
    -\zeta'_{1, d}(0) + \sum_{q \geq 2} (-1)^q q \zeta'_{q, d}(0).
$$
By virtue of Theorem~2.1,
we have a positive constant $C_1$ such that
$$
\sum_{q \geq 2} (-1)^q q \zeta'_{q, d}(0) \leq C_1 d^n \log(d).
$$
On the other hand, by Proposition 2.7.6 of \cite{Vo},
for some positive $C_2$,
we get
$$
    \zeta'_{1, d}(0) \geq -C_2 d^n \log(d).
$$
Thus we obtain our corollary.
\qed
\enddemo
\vfill\eject

\head 3. Hermitian modules over arithmetic curves
\endhead

In this section, we will discuss Hermitian modules over arithmetic curves.
Some results are found in \cite{St} and \cite{Gr}. For reader's convenience,
we will give however an explicit proof for each result.

Let $V$ be a $\CC$-vector space, $h_V$ a Hermitian metric on $V$ and
$W$ a subvector space of $V$. The metric $h_V$ induces a metric $h_W$ of
$W$, which is called {\it the submetric of $W$ induced by $h_V$}.
Let $W^{\perp}$ be the orthogonal complement of $W$.
Then the natural homomorphism $W^{\perp} \longrightarrow V/W$ is isomorphic.
Thus we have a metric $h_{V/W}$ of $V/W$. This metric is called
{\it the quotient metric of $V/W$ induced by $h_V$}.

\proclaim{Lemma 3.1}
With notation as above,
let $x_1, \ldots, x_s$ be elements of $W$ and $x_{s+1}, \ldots , x_n$
elements of $V$ such that
$\{ x_1, \ldots , x_s \}$ is a basis of $W$,
$\{ x_1, \ldots , x_n \}$ is a basis of $V$ and
$\{ \bar{x}_{s+1} , \ldots , \bar{x}_n \}$ is a basis of $V/W$,
where $\bar{x}_{s+1} , \ldots , \bar{x}_n$ are images of
$x_{s+1} , \ldots, x_{n}$ in $V/W$. Then we have
$$
\align
 \split
 \vmatrix
   h_{V}(x_1, x_1) & \cdots & h_{V}(x_1, x_n) \\
   \vdots                    &        & \vdots    \\
   h_{V}(x_n, x_1) & \cdots & h_{V}(x_n, x_n)
 \endvmatrix = &
 \vmatrix
   h_{W}(x_1, x_1) & \cdots & h_{W}(x_1, x_s) \\
   \vdots                    &        & \vdots    \\
   h_{W}(x_s, x_1) & \cdots & h_{W}(x_s, x_s)
 \endvmatrix \times \\
 &
 \vmatrix
   h_{V/W}(\bar{x}_{s+1}, \bar{x}_{s+1}) & \cdots &
   h_{V/W}(\bar{x}_{s+1}, \bar{x}_n)   \\
   \vdots                    &        & \vdots \\
   h_{V/W}(\bar{x}_n, \bar{x}_{s+1}) & \cdots &
   h_{V/W}(\bar{x}_n, \bar{x}_n)
 \endvmatrix
 \endsplit
\endalign
$$
\endproclaim

\demo{Proof}
Let $x_i = y_i + z_i$ be decompositions such that
$y_i \in W$ and $z_i \in W^{\perp}$.
Then it is easy to see that
$$
 \vmatrix
   h_{V/W}(\bar{x}_{s+1}, \bar{x}_{s+1}) & \cdots &
   h_{V/W}(\bar{x}_{s+1}, \bar{x}_n)   \\
   \vdots                    &        & \vdots \\
   h_{V/W}(\bar{x}_n, \bar{x}_{s+1}) & \cdots &
   h_{V/W}(\bar{x}_n, \bar{x}_n)
 \endvmatrix =
 \vmatrix
   h_{V}(z_{s+1}, z_{s+1}) & \cdots & h_{V}(z_{s+1}, z_n)\\
   \vdots                            &        & \vdots \\
   h_{V}(z_n, z_{s+1}) & \cdots & h_{V}(z_n, z_n)
 \endvmatrix
$$
and
$$
 (x_1, \ldots , x_n) = (x_1, \ldots , x_s, z_{s+1} , \ldots , z_n) U,
 \quad\hbox{where}\quad
 U = \pmatrix
  I_s & * \\
  0   & I_{n-s}
 \endpmatrix .
$$
Hence
$$
\align
 \split &
 \pmatrix
   h_{V}(x_1, x_1) & \cdots & h_{V}(x_1, x_n) \\
   \vdots                    &        & \vdots    \\
   h_{V}(x_n, x_1) & \cdots & h_{V}(x_n, x_n)
 \endpmatrix = \\
 &
 U \pmatrix
   h_{W}(x_1, x_1) & \cdots & h_{W}(x_1, x_s)  & 0 & \cdots & 0 \\
   \vdots          &        & \vdots           &\vdots&     &\vdots \\
   h_{W}(x_s, x_1) & \cdots & h_{W}(x_s, x_s)  & 0 & \cdots & 0 \\
   0 & \cdots & 0 & h_{V}(z_{s+1}, z_{s+1}) & \cdots & h_{V}(z_{s+1}, z_n)\\
   \vdots & &\vdots&\vdots                  &        & \vdots \\
   0 & \cdots & 0 & h_{V}(z_n, z_{s+1}) & \cdots & h_{V}(z_n, z_n)
 \endpmatrix
 {}^{t}\bar{U}.
 \endsplit
\endalign
$$
Thus we have our lemma.
\qed
\enddemo

\proclaim{Proposition 3.2}
Let $K$ be an algebraic number field and
$O_K$ the ring of integers of $K$.
Let $0 \to S \to E \to Q \to 0$ be an exact sequence of
$O_K$-modules and $h'$, $h$ and $h''$ Hermitian metrics of $S$, $E$ and
$Q$ respectively. If, for each infinite place $\sigma$ of $K$,
$h'_{\sigma}$ is the submetric of $h_{\sigma}$ and
$h''_{\sigma}$ is the quotient metric of $h_{\sigma}$, then
we have
$$
  \deg_{L^2}(E, h) = \deg_{L^2}(S, h') + \deg_{L^2}(Q, h'').
$$
\endproclaim

\demo{Proof}
This is an immediate consequence of Lemma~3.1.
\qed
\enddemo

\bigskip
Let $K$ be an algebraic number field and
$O_K$ the ring of integers of $K$.
Let $(E, h)$ be a Hermitian module over $O_K$.
We define an averaged $L^2$-degree $\mu_{L^2}(E, h)$ of $(E, h)$ as follows:
$$
    \mu_{L^2}(E, h) = \frac{\deg_{L^2}(E, h)}{\rank E}.
$$

\proclaim{Proposition 3.3}
With notation as above, we set $S_K = \Spec(O_K)$.
Let $h'$ be another Hermitian metrics of $E$.
Assume that, for each $\sigma \in (S_K)_{\sigma}$, there
is a constant $C_{\sigma}$ such that
$$
    h_{\sigma}(x, x) \leq C_{\sigma} h'_{\sigma}(x, x)
$$
for all $x \in E_{\sigma}$. Then we have
$$
\mu_{L^2}(E, h') \leq \mu_{L^2}(E, h) + \frac{1}{2}
\sum_{\sigma \in (S_K)_{\infty}}  \log C_{\sigma}
$$
\endproclaim

\demo{Proof}
Clearly it is sufficient to show the following lemma.
\qed
\enddemo

\proclaim{Lemma 3.4}
Let $V$ be a vector space over the complex number field $\CC$ and
$h$ and $h'$ Hermitian metrics on $V$ such that $h(x, x) \leq h'(x, x)$ for
all $x \in V$. Let $e_1, \ldots, e_n$ be a basis of $V$. Then we have
$\det (h(e_i, e_j)) \leq \det (h'(e_i, e_j))$. Moreover the equality
holds if and only if $h = h'$.
\endproclaim

\demo{Proof}
We set $H = (h(e_i, e_j))$ and $H' = (h'(e_i, e_j))$.
Let $U$ be an unitary matrix such that
$h'(Ue_i, Ue_j) = \lambda_i \delta_{ij}$.
Here we put $h_{ij} = h(Ue_i, Ue_j)$. Then
by virtue of Hadamard's inequality we have
$$
 \det H = \det {}^t U H \overline{U} \leq h_{11} \cdots h_{nn}
 \leq \lambda_1 \cdots \lambda_n = \det {}^t U H' \overline{U} = \det H'.
$$
Furthermore, if the equality holds, then $h_{ij} = 0$ for all $i \not= j$
and $h_{ii} = \lambda_i$ for all $i$
by the equality condition of Hadamard's inequality. Thus
${}^t U H \overline{U} = {}^t U H' \overline{U}$. Hence $h = h'$.
\qed
\enddemo

\bigskip
Let $k \subseteq K$ be an extension of algebraic number fields and
$O_k$ and $O_K$ the ring of integers of $k$ and $K$ respectively.
We set $S_k = \Spec(O_k)$ and $S_K = \Spec(O_K)$
and denote by $f$ the natural morphism $S_K \longrightarrow S_k$.
If for $\sigma \in \left( S_k \right)_{\infty}$ we set
$$
 f_{\infty}^{-1}(\sigma) = \{ \sigma' \in \left( S_K \right)_{\infty} \mid
                          \rest{\sigma'}{k} = \sigma \},
$$
we have a natural isomorphism
$$
   f_*(V) \emotimes{\sigma}{k} \CC \overset{\sim}\to\longrightarrow
   \bigoplus_{\sigma' \in f_{\infty}^{-1}(\sigma)}
   V \emotimes{\sigma'}{K} \CC,
$$
where $f_*(V) \emotimes{\sigma}{k} \CC$ is a tensor product
by the embedding $\sigma : k \hookrightarrow \CC$.
Thus we have a metric $h'(\sigma)$ on $f_*(V) \emotimes{\sigma}{k} \CC$
by the above isomorphism.
We denote this Hermitian module on $S_k$ by $(f_* V, f_* h)$ and call it
the push-forward of $(V, h)$. The Riemann-Roch theorem in this situation
asserts
$$
\deg_{L^2}(f_*V, f_*h) = \deg_{L^2}(V, h) + (\rank V)
\deg_{L^2}(f_* O_K, f_* \operatorname{can}_K).
$$

\proclaim{Proposition 3.5}
Let $(E, h)$ be a Hermitian module on $S_K$. Then the set
$$
 \{ \mu_{L^2}(F, h_F) \mid
 \hbox{$F$ is a sub-sheaf of $E$ and $h_F$ is the induced
                           metric by $h$.} \}
$$
is discrete subset of $\RR$ and bounded as above.
\endproclaim

\demo{Proof}
We set
$$
    M_l(E, h) =
    \left\{ \deg_{L^2}(F, h_F) \left|
    \matrix
    \hbox{$F$ is a sub-sheaf of $E$ with $\rank(F) = l$ and} \hfil \\
    \hfil \hbox{$h_F$ is the induced metric by $h$.} \hfill
    \endmatrix
    \right. \right\}
$$
Then we have
$$
\align
 \split
 &
 \{ \mu_{L^2}(F, h_F) \mid
    \hbox{$F$ is a sub-sheaf of $E$ and $h_F$ is the induced
          metric by $h$.} \} = \\
 & \qquad
 M_0(E, h) \cup M_1(E, h) \cup \cdots \cup \frac{1}{l} M_l(E, h)
  \cup \cdots \cup \frac{1}{\rank E} M_{\rank E}(E, h).
 \endsplit
\endalign
$$
Thus it is sufficient to see that $M_{l}(E, h)$ is discrete on $\RR$ and
bounded as above for each $l$.
Let $E_{tor}$ be the torsion part of $E$ and $\bar{E} = E/E_{tor}$.
Let $F$ be a sub-sheaf of $E$ with $\rank F = l$,
$F_{tor}$ the torsion part of $F$ and $\bar{F} = F/F_{tor}$.
Then,
$$
    \deg_{L^2}(F, h_F) = \log\#(F_{tor}) + \deg_{L^2}(\bar{F}, h_{\bar{F}}),
    \quad
    \bar{F} \subseteq \bar{E}
    \quad\hbox{and}\quad
    F_{tor} \subseteq E_{tor}.
$$
Thus we have
$$
   M_l(E, h) \subseteq M_l(\bar{E}, h) + \{ \log(i) \mid
   \hbox{$i \in \ZZ$ with $1 \leq i \leq \#(E_{tor})$} \}
$$
Therefore we may assume $E$ is torsion free.

If we denote by $f$ the natural morphism $S_K \longrightarrow S_{\QQ}$,
by Riemann-Roch theorem, we have
$$
\deg_{L^2}(f_* F, f_* h_F) = \deg_{L^2}(F, h_F) +
l \deg_{L^2}(f_* O_K, f_* \operatorname{can}_K)
$$
implies
$$
   M_l(E, h) \subseteq M_{ln}(f_*E, f_*h) -
                       l\deg_{L^2}(f_* O_K, f_*\operatorname{can}_K),
$$
where $n = [K : \QQ]$. Thus we may assume $K = \QQ$.
Since $M_l(E, h) \subseteq M_1(\bigwedge^{l} E, \bigwedge^{l} h)$,
we may furthermore assume $l = 1$.

We set $E = \ZZ x_1 \oplus \cdots \oplus \ZZ x_r$ and
$$
 H =
  \pmatrix
   h(x_1, x_1) & \cdots & h(x_1, x_r) \\
   \vdots      &        & \vdots      \\
   h(x_r, x_1) & \cdots & h(x_r, x_r)
  \endpmatrix.
$$
Since $H$ is a Hermitian matrix, there is an unitary matrix $U$ such that
$$
   H =
   {}^t U \diag(\lambda_1, \ldots, \lambda_r) \bar{U}
$$
and $0 < \lambda_1 \leq \cdots \leq \lambda_r$.
Let $F$ be a sub-module of $E$ of rank $1$.
We set $F = \ZZ (a_1 x_1 + \cdots + a_r x_r)$,
$a_i \in \ZZ$ $(1 \leq i \leq r)$ and
$(b_1, \ldots, b_r) = (a_1, \ldots, a_r){}^{t}U$.
Then
$$
\align
   \deg_{L^2}(F, h_F) & =
     - \frac{1}{2} \log h((a_1, \ldots, a_r), (a_1, \ldots, a_r)) \\
     & =
     - \frac{1}{2} \log( \lambda_1 |b_1|^2 + \cdots + \lambda_r |b_r|^2)
\endalign
$$
Thus since $|a_1|^2 + \cdots + |a_r|^2 = |b_1|^2 + \cdots + |b_r|^2$,
we have
$$
  -\frac{1}{2} \log( \lambda_r(|a_1|^2 + \cdots + |a_r|^2)) \leq
  \deg_{L^2}(F, h_F) \leq
  -\frac{1}{2} \log( \lambda_1(|a_1|^2 + \cdots + |a_r|^2))
$$
Thus ${\displaystyle \deg(F, h_F) \leq -\frac{1}{2} \log( \lambda_1)}$,
which shows
$M_1(E, h)$ is bounded as above.
Here we assume that $M_1(E, h)$ is not discrete in $\RR$.
Then there is a sequence $\{ F_n \}$ such that
$F_n$ are distinct rank $1$ sub-sheaves of $E$ and
${\displaystyle \lim_{n\to\infty} \deg_{L^2}(F_n, h_{F_n})}$ exits.
We set $F_n = \ZZ (a^{(n)}_1 x_1 + \cdots + a^{(n)}_r x_r)$.
Since $F_n$ are distinct, we have
$$
\lim_{n\to\infty} |a^{(n)}_1|^2 + \cdots + |a^{(n)}_r|^2
= \infty.
$$
Thus by the above inequality,
${\displaystyle \lim_{n\to\infty} \deg_{L^2}(F_n, h_{F_n}) = -\infty}$.
This is a contradiction.
\qed
\enddemo

\proclaim{Proposition 3.6}
Let $k \subseteq K$ be an extension of algebraic number fields and
$O_k$ and $O_K$ the ring of integers of $k$ and $K$ respectively.
We denote by $f$ the natural morphism
$\Spec(O_K) \longrightarrow \Spec(O_k)$.
For a Hermitian module $(L, h)$ on $\Spec(O_K)$ of rank $1$,
there is a constant $C$
such that
$$
   \mu_{L^2}(F, h_F) \leq C m
$$
for all $m \geq 1$ and all $O_k$-submodule $F$ of $f_* L^m$,
where $h_F$ is the induced metric by $f_* h^m$.
\endproclaim

\demo{Proof}
First we claim that we may assume $L = O_{K}$.
Let $x$ be a non-zero element of $L$. We define
a homomorphism $\phi_m : O_{K} \longrightarrow L^m$ by
$\phi_m(1) = x \otimes \cdots \otimes x$. Let $h_1$ be the induced metric
of $O_{K}$ by $\phi_1$. Then clearly the induced metric by $\phi_m$ is
equal to $h_1^m$.
We set $Q_m = \Coker (\phi_m : O_{K} \longrightarrow L^m)$.
We consider the exact sequence:
$$
  0 \to f_*O_{K} \to f_* L^m \to f_* Q_m \to 0.
$$
Let $F$ be a $O_k$-submodule of $f_* L^m$. Then,
$$
  \deg_{L^2}(F) \leq \deg_{L^2}(F \cap f_* O_{K}) + \log \#(Q_m).
$$
Thus we have the claim.

Let $F$ be a $O_k$-submodule of $f_* O_{K}$.
Let $h_F$ be the induced metric by $f_* h^m$ and $h'_F$ the induced metric
by $f_* \operatorname{can}_{K}$.
Then, it is easy to see that
$$
  \mu_{L^2}(F, h_F) = \mu_{L^2}(F, h'_F) -
                \frac{m}{2} \sum_{\sigma \in \left( S_K \right)_{\infty}}
                  \log h_{\sigma}(1, 1).
$$
Hence we have our proposition by Proposition~3.5.
\qed
\enddemo
\vfill\eject

\head 4. Asymptotic behavior of $L^2$-degree of submodules of $H^0(L^m)$
\endhead

A main purpose of this section is to prove the following theorem.

\proclaim{Theorem 4.1}
Let $K$ be an algebraic number field, $O_K$ the ring of integers of $K$
and $f : X \longrightarrow \Spec(O_K)$ an arithmetic variety over $O_K$.
Let $(L, h)$ be a Hermitian line bundle on $X$,
$(E, h_E)$ a Hermitian vector bundle on $X$, and
$h_m$ a Hermitian metric of $H^0(X, L^m \otimes E)$ induced by
$h^m \otimes h_E$.
Then there is a constant $C$ such that
$$
    \frac{\deg_{L^2}(F, h_F)}{\rank F} \leq C m \log m
$$
for all $m > 1$ and all submodules $F$ of $H^0(X, L^m \otimes E)$,
where $h_F$ is the Hermitian metric induced by $h_m$.
\endproclaim

First of all, we prepare the following Lemma.

\proclaim{Lemma 4.2}
Let $M$ be a compact K\"{a}hler manifold of dimension $d$,
$(L, h)$ a Hermitian line bundle
on $M$ and $\Phi$ a K\"{a}hler form on $M$.
Let $s \in H^0(M, L)$ be a smooth section and
$V$ the set of zero points of $s$.
Then there are constants $C$ and $C'$ such that
$$
   \int_M h^{m-1} (t , t) \Phi^d \leq
   C m^{6 d} \int_M h^{m}(s \otimes t, s \otimes t) \Phi^d
$$
for all $m \geq 1$ and $t \in H^0(M, L^{m-1})$ and that
$$
   \int_V h^m(t, t) \Phi^{d-1} \leq C' m^{2 d} \int_M h^m(t, t) \Phi^d
$$
for all $m \geq 1$ and $t \in H^0(M, L^{m})$.
\endproclaim

\demo{Proof}
Here we claim the following sublemma, which can be proved by the same
idea (due to M. Gromov) as in the proof of Proposition 3 of \cite{GS88}.

\proclaim{Sublemma 4.2.1}
For a real number $p \geq 1$, there is a constant $C$ such that
$$
     || t ||_{\sup} \leq C m^{2d/p} || t ||_{L^p}
$$
for all $m \geq 1$ and $t \in H^0(M, L^m)$.
\endproclaim

\demo{Proof}
Let $x$ be a point of $M$, $U_x$ an open neighborhood of $x$,
$$
  U_x \overset\sim\to\longrightarrow V_x \subseteq \CC^d
$$
a local chart, and $u_x$ a local basis of $L$.
{}From now, we identify $U_x$ with $V_x$.
We set $q_x(z) = h(u_x, u_x)$ and
$$
     r(a, b) = |a_1 - b_1| + \cdots + |a_d - b_d|
$$
for $a = (a_1, \ldots, a_d), b = (b_1, \ldots, b_d) \in U_x$.
Moreover we set
$$
  B(a, R) = \{ (z_1, \ldots, z_d) \in U_x \mid
  \hbox{$|z_i - a_i| \leq R$ for all $i$} \}
$$
for a point $a = (a_1, \ldots, a_d) \in U_x$ and $R > 0$.
Shrinking $U_x$ if necessary, we may assume that
there is a constant $K_x$ such that
$$
     |q_x(a) - q_x(b)| \leq K_x r(a, b)
$$
for all $a, b \in U_x$.
Thus there are an open neighborhood $W_x \subseteq U_x$
of $x$ and a constant $r_x$ with the following properties:
\roster
\item "(1)" $B(a, r_x) \subseteq U_x$ for all $a \in W_x$.

\item "(2)" $q_x(z) \geq q_x(a) - K_x r(a, z) > 0$ for $z \in B(a, r_x)$ and
$a \in W_x$.
\endroster
Since $M = \bigcup_{x \in M} W_x$ and $M$ is compact,
there are finitely many $x_1, \ldots, x_n$ such that
$M = \bigcup_{i=1}^{n} W_{x_i}$.
For simplicity, we set $U_i = U_{x_i}$, $W_i = W_{x_i}$, $u_i = u_{x_i}$,
$q_i = q_{x_i}$, $K_i = K_{x_i}$ and $r_i = r_{x_i}$.
Let $t$ be a section of $H^0(M, L^m)$.
Since $M$ is compact, we can take $a \in M$ with
$||t||_{\sup} = \sqrt{h^m(t, t)(a)}$.
Let $a \in W_i$ ($1 \leq i \leq n$).
Using local basis $u_i$, we write $t = f(z)u_i^m$.
Then,
$$
    h^m(t, t)^{p/2} = |f(z)|^p q_i(z)^{pm/2}.
$$
We take constants $g_i$, $G_i$ and $C_i$ such that
$g_i \leq q_i(z) \leq G_i$ for all $z \in U_i$ and
$\Phi^d \geq C_i dx_1 dy_1 \cdots dx_d dy_d$, where
$z_j = x_j + \sqrt{-1} y_j$ for $1 \leq j \leq d$.
Let $l_m$ be the integral part of $pm/2$ and $l'_m = pm/2 - l_m$.
Since $|f|^p$ is subharmonic, we have
$$
\align
||t||_{L^p}^p & = \int_M h^m(t, t)^{p/2} \Phi^d \\
& \geq C_i \int_{B(a, r_i)}
|f(z)|^p q_i(z)^{pm/2} dx_1 dy_1 \cdots dx_d dy_d \\
& \geq C_i g_i^{l'_m}
\int_{B(a, r_i)} |f(z)|^p q_i(z)^{l_m} dx_1 dy_1 \cdots dx_d dy_d \\
& \geq C_i g_i^{l'_m} \int_{B(a, r_i)} |f(z)|^p (q_i(a) - K_i r(a, z))^{l_m}
dx_1 dy_1 \cdots dx_d dy_d \\
& \geq (2 \pi)^d C_i g_i^{l'_m} |f(a)|^p \times \\
& \phantom{\geq} \qquad \int_0^{r_i} \dotsi \int_0^{r_i}
R_1 \cdots R_d (q_i(a) - K_i(R_1 + \cdots + R_d))^{l_m}
dR_1 \cdots dR_d.
\endalign
$$
On the other hand, by an easy calculation, we can see
$$
\int_0^r x ( A - B x)^l dx \geq \frac{r^2}{(l + 1)(l + 2)} A^l,
$$
where $A > 0$, $B > 0$, $r > 0$, $A - B r > 0$ and $l$ is a non-negative
integer. Hence we have
$$
\align
||t||_{L^p}^p & \geq \frac{(2 \pi)^d C_i g_i^{l'_m} r_i^{2d}}
                          {(l_m + 1)^d(l_m + 2)^d} |f(a)|^p q_i(a)^{l_m} \\
& \geq \frac{(2 \pi)^d C_i r_i^{2d}}
          {(l_m + 1)^d(l_m + 2)^d}
    \left( \frac{g_i}{G_i} \right)^{l'_m}
    |f(a)|^p q_i(a)^{pm/2} \\
& \geq \frac{(2 \pi)^d C_i r_i^{2d}}
            {(pm/2 + 1)^d(pm/2 + 2)^d}
    \left( \frac{g_i}{G_i} \right) ||t||_{\sup}^p.
\endalign
$$
Here we take a constant $C$ such that
$$
\frac{(2 \pi)^d C_i r_i^{2d}}{(pm/2 + 1)^d(pm/2 + 2)^d}
\left( \frac{g_i}{G_i} \right) \geq \frac{C}{m^{2d}},
$$
for all $m > 0$ and $1 \leq i \leq n$.
Then
$$
||t||_{L^p}^p \geq \frac{C}{m^{2d}} ||t||_{\sup}^p.
$$
Thus we have our sublemma.
\qed
\enddemo

\bigskip
To prove the first part of Lemma~4.2,
it is sufficient to see the following claim.

\proclaim{Claim 4.2.2}
There are constants $C_1$ and $C_2$ such that
$$
     || t ||_{L^1} \leq C_1 m^{d} || s \otimes t ||_{L^2}
$$
for all $m \geq 1$ and $t \in H^0(M, L^{m-1})$ and that
$$
      || t ||_{L^2} \leq C_2 m^{2d} || t ||_{L^1}
$$
for all $m \geq 1$ and $t \in H^0(M, L^{m-1})$.
\endproclaim

Since $h(s, s)^{-1/2}$ is integrable, we have
$$
\align
\int_M h^{m-1}(t, t)^{1/2} \Phi^d & =
\int_M h(s, s)^{-1/2} h^m(s \otimes t, s \otimes t)^{1/2} \Phi^d \\
& \leq \left( \int_M h(s, s)^{-1/2} \Phi^d \right)
       || s \otimes t ||_{\sup}.
\endalign
$$
Thus by Sublemma~4.2.1, we have first inequality of Claim~4.2.2.

Clearly there is a constant $C_3$ such that
$$
      || t ||_{L^2} \leq C_3  || t ||_{\sup}
$$
Therefore Sublemma~4.2.1 implies second inequality of Claim~4.2.2.

\bigskip
Next we consider the second assertion of Lemma~4.2.
Let $v$ be the volume of $V$.
Using Sublemma~4.2.1, for some constant $C_4$, we have
$$
\align
   \int_V h^m(t, t) \Phi^{d-1} & \leq v ||\rest{t}{V}||_{\sup}^2 \\
   & \leq v || t ||_{\sup}^2 \\
   & \leq v C_4 m^{2 d} \int_M h^m(t, t) \Phi^d.
\endalign
$$
Thus we get Lemma~4.2.
\hfill\qed
\enddemo

\bigskip
Now let us go back to the proof of Theorem~4.1.
To get Theorem~4.1, clearly we may assume that $X$ is normal.

\subhead 4.3 Reduction 1 \endsubhead
We may assume that $E$ is a line bundle.

Let $g : Y = \PP(E) \longrightarrow X$ be the projective bundle of $E$
and $\OO(1)$ the tautological line bundle of $E$.
We can define the quotient metric $h_{\OO(1)}$ on $\OO(1)$ by
the natural surjective homomorphism $g^*(E) \to \OO(1)$.
We consider an isomorphism:
$$
\alpha_m : H^0(X, L^m \otimes E) \longrightarrow
           H^0(Y, g^*(L)^m \otimes \OO(1)).
$$
Let $F$ be a submodule of $H^0(X, L^m \otimes E)$.
Let $h_F$ be a Hermitian metric of $F$ induced by $h^m \otimes h_E$ and
$h_{\alpha_m(F)}$ a Hermitian metric of $\alpha_m(F)$ induced by
$g^*(h)^m \otimes h_{\OO(1)}$. Then, by Proposition~3.3, we have
$$
    \mu(F, h_F) \leq \mu(\alpha_m(F), h_{\alpha_m(F)}).
$$
Thus we get Reduction 1.

\subhead 4.4 Reduction 2 \endsubhead
We may assume that $H^0(X, E) \not= 0$.

Let $(E', h_{E'})$ be a Hermitian line bundle on $X$ such that
$H^0(X, E') \not= 0$ and $H^0(X, E \otimes E') \not= 0$.
Let $s_1$ be a non-zero section of $H^0(X, E')$.
We consider an injective homomorphism:
$$
\beta_m : H^0(X, L^m \otimes E) \overset{\otimes s_1}\to\longrightarrow
          H^0(X, L^m \otimes E \otimes E').
$$
We take a constant $B_1$ such that
$||(s_1)_{\sigma}||_{\sup} \leq B_1$
for all $\sigma \in K_{\infty}$.
Then we have
$$
       ||{\beta_m(t)}_{\sigma}||_{L^2} \leq B_1 ||t_{\sigma}||_{L^2}
$$
for all $t \in H^0(X, L^m \otimes E)$ and all $\sigma \in K_{\infty}$.
Thus, by Proposition~3.3, for a submodule $F$ of
$H^0(X, L^m \otimes E)$, we obtain
$$
\mu(F, h_F) \leq \mu(\beta_m(F), h_{\beta_m(F)}) +
[K : \QQ] \log B_1.
$$
Therefore, we get Reduction 2.

\subhead 4.5 Reduction 3 \endsubhead
We may assume that $E = (\OO_X, \operatorname{can}_K)$.

Let $s_2$ be a non-zero section of $H^0(X, E)$.
We consider an injective homomorphism:
$$
\gamma_m : H^0(X, L^m \otimes E) \overset{\otimes s_2^{m-1}}\to\longrightarrow
          H^0(X, (L \otimes E)^m).
$$
We take a constant $B_2$ such that
$||(s_2)_{\sigma}||_{\sup} \leq B_2$
for all $\sigma \in K_{\infty}$.
Then, by the same way as in Reduction 2, for a submodule $F$ of
$H^0(X, L^m \otimes E)$, we have
$$
\mu(F, h_F) \leq \mu(\gamma_m(F), h_{\gamma_m(F)}) +
(m-1)[K : \QQ] \log B_2.
$$
Therefore, we get Reduction 3.

\subhead 4.6 Reduction 4 \endsubhead
We may assume that $L$ is very ample.

Let $(L', h_{L'})$ be a Hermitian line bundle on $X$ such that
$H^0(X, L') \not= 0$ and $L \otimes L'$ is very ample.
Let $s_3$ be a non-zero section of $H^0(X, L')$.
We consider an injective homomorphism:
$$
\delta_m : H^0(X, L^m) \overset{\otimes s_3^{m}}\to\longrightarrow
           H^0(X, (L \otimes L')^m).
$$
We take a constant $B_3$ such that
$||(s_3)_{\sigma}||_{\sup} \leq B_3$
for all $\sigma \in K_{\infty}$.
Then, by the same way as in Reduction 2, for a submodule $F$ of
$H^0(X, L^m \otimes E)$, we have
$$
\mu(F, h_F) \leq \mu(\delta_m(F), h_{\delta_m(F)}) +
m[K : \QQ] \log B_3.
$$
Therefore, we get Reduction 4.

\subhead 4.6 \endsubhead
Let us go to the main part of the proof of Theorem~4.1.
Gathering Reduction 1--4, we may assume that
$(E, h) = (\OO_X, \operatorname{can}_K)$ and $L$ is very ample.
We use induction on $\dim f$. If $\dim f = 0$,
this theorem holds by Proposition~3.6.
Here we remark the following.

\definition{Remark 4.7.1}
For a base extension
$\Spec(O_{K'}) \longrightarrow \Spec(O_K)$ and
$O_K$-submodule $F \subseteq H^0(X, L^m)$,
we get
$$
F \otimes O_{K'} \subseteq
H^0(X \otimes O_{K'}, (L \otimes O_{K'})^m)
$$
and
$$
\mu_{L^2}(F \otimes O_{K'}, h_F \otimes O_{K'}) =
[K' : K] \mu_{L^2}(F, h_F).
$$
\enddefinition

For $s \in H^0(X, L)$, we denote by $\operatorname{div}(s)$
the set of zero of $s$.
Since $L$ is very ample, by the above Remark~4.7.1,
considering a base change of $f : X \longrightarrow \Spec(O_K)$
if necessarily, we may assume that there is a section $s \in H^0(X, L)$
such that
$\operatorname{div}(s) \longrightarrow \Spec(O_K)$ is generically smooth.
Let $V' = \operatorname{div}(s)$, $V$ the horizontal component of $V'$,
$I_V$ the defining ideal of $V$.
We give a Hermitian metric $h_{I_V}$ of $I_V$ by $h^{-1}$.
Let
$$
X \overset{f'}\to\longrightarrow \Spec(O_{K'})
\overset{\pi}\to\longrightarrow \Spec(O_K)$$
be the Stein factorization of $X \overset{f}\to\longrightarrow \Spec(O_K)$.
Then $f'_*(I_V \otimes \OO_X(V'))$ is a torsion
free sheaf of rank $1$ and the natural homomorphism
$$
I_V \otimes \OO_X(V') \longrightarrow {f'}^*(f'_*(I_V \otimes \OO_X(V')))
$$
is injective
because $I_V \otimes \OO_X(V')$ is isomorphic to $\OO_{X_{K'}}$ on the
generic fiber $X_{K'}$ of $f'$.
Since $f'_*(I_V \otimes \OO_X(V'))$ is of rank 1, there is a factorial ideal
$A$ of $O_K$ such that $f'_*(I_V \otimes \OO_X(V')) \subseteq \pi^*(A)$.
In particular, we have $I_V \otimes \OO_X(V') \subseteq f^*(A)$.
Here we give a trivial Hermitian metric $h_A$ to $A$. Then the inclusion
$I_V \otimes \OO_X(V') \subseteq f^*(A)$ is isometry.

Let $F$ be a submodule of $H^0(X, L^m)$ and $h_F$ the induced metric by $h_m$.
We consider the exact sequence:
$$
0 \to H^0(X, L^{m} \otimes I_V) \to H^0(X, L^m)
  \to H^0(V, L^m).
$$
We set
$T = F \cap H^0(X, L^{m} \otimes I_V)$ and
$Q = F/F \cap H^0(X, L^{m} \otimes I_V)$.
Let $h_T$ be the submetric of $h_F$ and $h_Q$ the quotient metric of $h_F$.
Then by Proposition~3.2,
we have
$$
    \mu_{L^2}(F, h_F) = \frac{t}{t + q}\mu_{L^2}(T, h_T) +
                        \frac{q}{t + q}\mu_{L^2}(Q, h_Q),
\tag 4.7.2
$$
where $t = \rank T$ and $q =  \rank Q$.
Let $h'_T$ be the submetric induced by $h^m \otimes h_{I_V} = h^{m-1}$.
Since $Q$ is a submodule of $H^0(V, L^{m})$,
we get a metric $h'_Q$ induced by $\left( \rest{h}{V} \right)_m$.
If we set $d = \dim f$, by Proposition~3.3 and
Lemma~4.2, there are constants
$C_1$ and $C_2$ such that
$$
  \mu_{L^2}(T, h_T) \leq \mu_{L^2}(T, h'_T) + \frac{\log C_1 + 6 d \log m}{2}
  [K : \QQ]
\tag 4.7.3
$$
and
$$
  \mu_{L^2}(Q, h_Q) \leq \mu_{L^2}(Q, h'_Q) + \frac{\log C_2 + 2 d \log m}{2}
  [K : \QQ],
\tag 4.7.4
$$
for all $m \geq 1$.
By hypothesis of induction, there is a constant $C_3$ such that
$$
    \mu_{L^2}(N, h_N) \leq C_3 m \log m
\tag 4.7.5
$$
for all $m > 1$ and all submodule $N$ of $H^0(V, L^m)$.
On the other hand, due to Proposition~3.5,
we can take a constant $C_4$ such that
$$
      \mu_{L^2}(P, h_P) \leq C_4
\tag 4.7.6
$$
for all submodule $P$ of $H^0(X, L^2)$.
Moreover
$$
\mu_{L^2}(T, h'_T) =
\mu_{L^2}(T \otimes A^{-1}, h'_T \otimes h_A^{-1}) + \mu_{L^2}(A, h_A)
$$
and $T \otimes A^{-1} \subseteq H^0(X, L^{m-1})$.
We choose $C$ satisfying
$$
 \cases
      {\displaystyle \mu_{L^2}(A, h_A) +
                     \frac{\log C_1 + 6 d \log m}{2}[K : \QQ]
                     \leq C \log m} & \\
      {\displaystyle C_3 m \log m +
                     \frac{\log C_2 + 2 d \log m}{2} [K : \QQ]
                     \leq C m \log m} & \\
      C_4 \leq C (2 \log 2) &
 \endcases
\tag 4.7.7
$$
for all $m > 0$.
Here we claim:

\proclaim{Claim 4.7.8}
$ \mu_{L^2}(F, h_F) \leq C m \log m$ for all $m > 1$ and
all submodules $F$ of $H^0(X, L^m)$.
\endproclaim

We prove this claim by induction on $m$.
By hypothesis of induction on $m$, we have
$$
    \mu_{L^2}(T \otimes A^{-1}, h'_T \otimes h_A^{-1})
    \leq C (m-1) \log(m-1).
$$
Thus by (4.7.3) and (4.7.7), we get
$$
    \mu_{L^2}(T, h_T) \leq C m \log m.
$$
On the other hand, (4.7.4), (4.7.5) and (4.7.7)
imply
$$
     \mu_{L^2}(Q, h_Q) \leq C m \log m.
$$
Hence by (4.7.2), we obtain
$$
  \mu_{L^2}(F, h_F) \leq \frac{t}{t + q} C m \log m +
                         \frac{q}{t + q} C m \log m
                    = C m \log m.
\eqno{\qed}
$$
\vfill\eject

\head 5. Vanishing of a certain Bott-Chern secondary class
\endhead

This section is devoted to prove the following lemma.

\proclaim{Lemma 5.1}
Let $C$ be a compact Riemann surface, $E$ a vector bundle of rank $r$ on $C$
and $h$ a projectively flat metric of $E$.
Let $f : Y = \PP(E) \longrightarrow C$ be the projective
bundle of $E$ and $\OO_Y(1)$ the tautological line bundle.
Let $\Cal{E} : 0 \to F \to f^*E \to \OO_Y(1) \to 0$ be the canonical
exact sequence.
We give the canonical Hermitian metrics on
$F$, $f^*E$ and $\OO_Y(1)$ induced by the Hermitian metric $h$ of $E$.
Then we have
$$
\int_{Y} \left\{
\left(
c_1(\OO_Y(1)) - \frac{1}{r} f^* (c_1(E))
\right)
\sum_{i=1}^r (-1)^i \tilde{c}_i(\Cal{E}) c_1(\OO_Y(1))^{r-i} \right\}
= 0.
$$
\endproclaim

\demo{Proof}
First we notice that we may assume that $\det E$ is divisible
by $r$ in $\Pic(C)$ considering a finite covering of $C$.
We set
$$
  \Phi(\Cal{E}) =
     \sum_{i=1}^r (-1)^i \tilde{c}_i(\Cal{E}) c_1(\OO_Y(1))^{r-i}.
$$
Let $L$ be a Hermitian line bundle on $C$.
We will see $\Phi(\Cal{E} \otimes L) = \Phi(\Cal{E})$.
Using the formula of (1.3.3.2) in \cite{GS90b}, we have
$$
\align
\Phi(\Cal{E} \otimes L) & =
    \sum_{i=1}^r (-1)^i \tilde{c}_i(\Cal{E} \otimes L)
                 \left( c_1(\OO_Y(1)) + f^*(c_1(L)) \right)^{r-i} \\
& = \sum_{i=1}^r (-1)^i
    \left(
    \tilde{c}_i(\Cal{E}) + (r-i+1) \tilde{c}_{i-1}(\Cal{E})f^*(c_1(L))
    \right) \times \\
&   \phantom{= \sum_{i=1}^r (-1)^i}
    \left(
    c_1(\OO_Y(1))^{r-i} + (r-i)c_1(\OO_Y(1))^{r-i-1}f^*(c_1(L))
    \right) \\
& = \Phi(\Cal{E}) + \sum_{i=1}^r (-1)^i
       (r-i)\tilde{c}_{i}(\Cal{E})c_1(\OO_Y(1))^{r-i-1}f^*(c_1(L)) + \\
&   \phantom{= \Phi(\bar{E}) +}
    \sum_{i=1}^r (-1)^i
    (r-i+1) \tilde{c}_{i-1}(\Cal{E})c_1(\OO_Y(1))^{r-i}f^*(c_1(L)) \\
& = \Phi(\Cal{E}).
\endalign
$$
Since $\det E$ is divisible by $r$ in $\Pic(C)$, we can take
a Hermitian line bundle $(L, h_L)$ such that
$(L, h_L)^{\otimes r} = (\det E, \det h)$. We set
$E' = E \otimes L^{-1}$. Then $c_1(E') = 0$.
These observations show us that we may assume that $c_1(E) = 0$.

Thus, to get our assertion, it is sufficient to see that
$$
 \int_Y
 \tilde{c}_{i}(\Cal{E})c_1(\OO_Y(1))^{r-i+1} = 0
$$
for $i = 1, \ldots, r$.
We denote the homomorphisms $F \longrightarrow f^*E$ and
$f^* E \longrightarrow \OO_Y(1)$
by $\iota$ and $\tau$ respectively.
Let $p : Y \times \PP^1 \longrightarrow Y$ and
$q : Y \times \PP^1 \longrightarrow \PP^1$
be the natural projections.
Let $\sigma \in H^0(\PP^1, \OO_{\PP^1}(1))$
be a non-zero section such that $\sigma(\infty) = 0$.
Let $h_{\OO_{\PP^1}(1)}$ be a metric of $\OO_{\PP^1}(1)$ such that
$\rest{h_{\OO_{\PP^1}(1)}(\sigma, \sigma)}{0} = 1$.
We set
$$
  \tilde{E} = \Coker \left(
        p^* \iota \oplus (\operatorname{id}_{p^* F} \otimes \sigma) :
             p^* F \longrightarrow p^*f^* E \oplus (p^* F \otimes q^*
\OO_{\PP^1}(1))
                     \right).
$$
Let $F^{\perp}$ be the orthogonal complement of $F$ in $f^*E$.
Using the natural bijective homomorphism
$$
   p^* F^{\perp} \oplus (p^* F \otimes q^* \OO_{\PP^1}(1))
   \longrightarrow \tilde{E},
$$
we give a metric $\tilde{h}$ on $\tilde{E}$.
Let $x, x' \in p^* f^* E$ and
$$
    x = x_1 + x_2 \ (x_1 \in p^* F, x_2 \in p^* F^{\perp}), \qquad
    x' = x'_1 + x'_2 \ (x'_1 \in p^* F, x'_2 \in p^* F^{\perp})
$$
the orthonormal decompositions. Then we have
$$
   \tilde{h}(x, x') - h_{E}(x, x') = (||\sigma||^2 - 1)h_{E}(x_1, x'_1).
$$
These observation shows us that
$\rest{(\tilde{E}, \tilde{h})}{Y \times \{ 0 \}}$ is isometric to
$(f^*E, h_E)$. On the other hand, clearly
$\rest{(\tilde{E}, \tilde{h})}{Y \times \{ \infty \}}$ is isometric to
$(F, h_F) \oplus (\OO_Y(1), h_{\OO_Y(1)})$.
Hence by the construction of $\tilde{c}_i(\Cal{E})$
(cf. 1.3)
$$
    \tilde{c}_i(\Cal{E}) = \int_{\PP^1} c_i(\tilde{E}, \tilde{h})
                                          \log |z|^2,
$$
where $z$ is an inhomogeneous coordinate of $\PP^1$.
Thus using Fubini's theorem, we have
$$
\align
 \int_Y
 \tilde{c}_{i}(\Cal{E})c_1(\OO_Y(1))^{r-i+1} & =
 \int_Y \left\{
    c_1(\OO_Y(1))^{r-i+1} \int_{\PP^1} c_i(\tilde{E}, \tilde{h}) \log |z|^2
          \right\} \\
 & = \int_{Y \times \PP^1}
      p^* c_1(\OO_Y(1))^{r-i+1} c_i(\tilde{E}, \tilde{h}) \log |z|^2 \\
 & = \int_{\PP^1} \left\{
      \log |z|^2 \int_Y p^* c_1(\OO_Y(1))^{r-i+1} c_i(\tilde{E}, \tilde{h})
                    \right\}.
\endalign
$$
Hence it is sufficient to see that
$$
  \int_Y p^* c_1(\OO_Y(1))^{r-i+1} c_i(\tilde{E}, \tilde{h}) = 0
$$
for $z \not= \infty$.
Since $E$ is flat, we can take a local flat frame
$s_0, \ldots, s_{r-1}$ of $E$
such that $h_{E}(s_i, s_j) = \delta_{ij}$.
Let $X_0, \ldots, X_{r-1}$
be a homogeneous coordinate of $Y = \PP(E)$ corresponding to
$s_0, \ldots, s_{r-1}$.
On $X_0 \not= 0$, we set $z_i = X_i/X_0$ and
$$
       g = \frac{s_0 + \bar{z}_1 s_1 + \cdots + \bar{z}_{r-1} s_{r-1}}
                {1 + |z_1|^2 + \cdots + |z_{r-1}|^2}.
$$
Then $g \in F^{\perp}$ and $\tau(g) = \tau(s_0)$. Since
$$
   s_0 = (s_0 - g) + g \quad\hbox{and}\quad
   s_i = (s_i - z_i g) + z_i g \ (i=1, \cdots, r-1)
$$
are orthogonal decompositions of $s_0$ and $s_i$ respectively,
we have
$$
\multline
 (\tilde{h}(s_i, s_j)) =
  ||\sigma||^2 (\delta_{ij}) + \\
     \frac{1 - ||\sigma||^2}{1 + |z_1|^2 + \cdots + |z_{r-1}|^2}
     \pmatrix
      1   & \bar{z}_1    & \bar{z}_2    & \cdots & \bar{z}_{r-1} \\
      z_1 & |z_1|^2      & z_1\bar{z}_2 & \cdots & z_1\bar{z}_{r-1} \\
      z_2 & z_2\bar{z}_1 & |z_2|^2      & \ddots & z_2\bar{z}_{r-1} \\
      \vdots & \vdots    & \vdots       & \cdots & \vdots   \\
      z_{r-1} & z_{r-1}\bar{z}_1 & z_{r-1}\bar{z}_2 & \cdots & |z_{r-1}|^2 \\
     \endpmatrix.
\endmultline
$$
Thus we get
$$
c_i(\tilde{E}, \tilde{h}) \in
   A^{0, 0} [dz, dz_1 , \ldots, dz_{r-1},
        d\bar{z}, d\bar{z}_1, \ldots, d\bar{z}_{r-1}].
$$
On the other hand, since
$$
 h_{E}(g, g) =  \frac{1}{1 + |z_1|^2 + \cdots + |z_{r-1}|^2},
$$
we obtain
$$
c_1(\OO_Y(1)) \in
   A^{0, 0} [dz_1 , \ldots, dz_{r-1}, d\bar{z}_1, \ldots, d\bar{z}_{r-1}].
$$
Therefore since
$p^* c_1(\OO_Y(1))^{r-i+1} c_i(\tilde{E}, \tilde{h})$ is
$(r+1, r+1)$-form and lies in
$$
   A^{0, 0} [dz, dz_1 , \ldots, dz_{r-1},
        d\bar{z}, d\bar{z}_1, \ldots, d\bar{z}_{r-1}],
$$
we obtain
$$
  p^* c_1(\OO_Y(1))^{r-i+1} c_i(\tilde{E}, \tilde{h}) = 0
$$
on the chart $X_0 \not= 0$. Hence we have our lemma.
\qed
\enddemo
\vfill\eject

\head 6. Donaldson's Lagrangian and arithmetic second Chern class
\endhead

Let $M$ be an $n$-dimensional complex manifold, and
$E$ a vector bundle of rank $r$ on $M$.
Let $h$ and $k$ be Hermitian metrics of $E$, and
$\{ h_t \}_{0 \leq t \leq 1}$ a $C^{\infty}$-deformation of Hermitian metrics
of $E$ such that $h_0 = h$ and $h_1 = k$.
Let $R_t$ be the curvature of $h_t$. We set
$$
  \align
   Q_1(E, h, k) & = \log(\det(h)/\det(k)), \\
   Q_2(E, h, k) & = \sqrt{-1} \int_0^1
   \tr(h_t^{-1} \cdot (\partial_t h_t) \cdot R_t) dt, \\
   \sch_2(E, h, k) & = \sch_2( 0 \to (E, k) \to (E, h) \to (0, 0) \to 0).
  \endalign
$$
By Lemma~(3.6) in Chapter~VI of \cite{Ko}, $Q_2(E, h, k)$ is uniquely
determined by $h$ and $k$ up to
$\partial(A^{0,1}) + \overline{\partial}(A^{1,0})$.
Then, by comparing $(1, 1)$-part of the formula in
Corollary 1.30, ii) of \cite{BGS}, we have the following lemma.
(Note that
$\sch_2(E, h, k) = \sch_2( 0 \to (E, h) \to (E, k) \to (0, 0) \to 0)$
in the sense of \cite{BGS}.)

\proclaim{Lemma 6.1}
${\displaystyle \sch_2(E, h, k) = -\frac{1}{2 \pi} Q_2(E, h, k)}$ modulo
$\partial(A^{0,1}) + \overline{\partial}(A^{1,0})$.
\endproclaim

Here we assume that $M$ is compact and K\"{a}hler. Let $\Phi$ be a
K\"{a}hler form of $M$.
The Donaldson's Lagrangian $DL(E, h, k ; \Phi)$ is
defined by
$$
DL(E, h, k ; \Phi) = \int_M (Q_2(E, h, k) - \frac{c}{n} Q_1(E, h, k) \Phi)
\wedge \frac{\Phi^{n-1}}{(n-1)!},
$$
where
$$
c = \frac{{\displaystyle 2 n \pi \int_M c_1(E) \wedge \Phi^{n-1}}}
         {{\displaystyle r \int_M \Phi^n}}.
$$
We fix a Hermitian metric $k$ of $E$.
Then, by Proposition~(3.37) in Chapter~VI of \cite{Ko},
the following are equivalent.
\roster
\item "(i)" $DL(E, h_0, k ; \Phi)$ gives the absolute minimal value of
$$\{ DL(E, h, k ; \Phi) \mid \hbox{$h$ is a Hermitian metric of $E$} \}.$$

\item "(ii)" $h_0$ is Einstein-Hermitian.
\endroster

\proclaim{Lemma 6.2}
Let $K$ be an algebraic number field and $O_K$ the ring of integers.
Let $f : X \longrightarrow \Spec(\OO_K)$ be an arithmetic variety with
$\dim X = d \geq 2$, and $(H, h_H)$ a Hermitian line bundle on $X$
such that, for each $\sigma \in K_{\infty}$,
$c_1(H_{\sigma}, h_{H_{\sigma}})$ gives a K\"{a}hler form $\Phi_{\sigma}$ on
an infinite fiber $X_{\sigma}$.
Let $E$ be a vector bundle on $X$, and
$h$ and $h'$ Hermitian metrics of $E$.
If $\det(h) = \det(h')$, then we have
$$
(\achern{2}{E, h} - \achern{2}{E, h'}) \cdot \achern{1}{H, h_H}^{d-2}
= \frac{(d-2)!}{2\pi} \sum_{\sigma \in K_{\infty}}
DL(E_{\sigma}, h_{\sigma}, h'_{\sigma} ; \Phi_{\sigma}).
$$
\endproclaim

\demo{Proof}
Since $\det(h) = \det(h')$, we have $\achern{1}{E, h} = \achern{1}{E, h'}$.
Thus we get
$$
\achern{2}{E, h} - \achern{2}{E, h'}
= -(\ach_2(E, h) - \ach_2(E, h')).
$$
Therefore, by Lemma~6.1, we obtain
$$
\achern{2}{E, h} - \achern{2}{E, h'} =
\frac{1}{2\pi}\sum_{\sigma \in K_{\infty}}
Q_2(E_{\sigma}, h_{\sigma}, h'_{\sigma}).
$$
Hence we have our formula because
$Q_1(E_{\sigma}, h_{\sigma}, h'_{\sigma}) = 0$ for all $\sigma \in K_{\infty}$.
\qed
\enddemo

\proclaim{Theorem 6.3}
Let $K$ be an algebraic number field and $O_K$ the ring of integers.
Let $f : X \longrightarrow \Spec(O_K)$ be an arithmetic variety with
$\dim X = d \geq 2$, and $(H, h_H)$ a Hermitian line bundle on $X$
such that, for each $\sigma \in K_{\infty}$,
$c_1(H_{\sigma}, h_{H_{\sigma}})$ gives a K\"{a}hler form $\Phi_{\sigma}$ on
an infinite fiber $X_{\sigma}$.
Let $E$ be a vector bundle of rank $r$ on $X$.
For a Hermitian metric $h$ of $E$, we set
$$ \Delta(E, h) = \left(
\achern{2}{E, h} - \frac{r - 1}{2r}\achern{1}{E, h}^2 \right)\cdot
\achern{1}{H, h_H}^{d-2}.$$
If $E$ is $\Phi_{\sigma}$-poly-stable on each infinite fiber $X_{\sigma}$,
then we have;
\roster
\item "(1)" the set
$\Delta = \{ \Delta(E, h) \mid \hbox{$h$ is a Hermitian metric of $E$} \}$
has the absolute minimal value.

\item "(2)" $\Delta(E, h_0)$ attaches the minimal value of $\Delta$
if and only if
$h_0$ is weakly Einstein-Hermitian on each infinite fiber.
\endroster
\endproclaim

\demo{Proof}
Since $E$ is poly-stable on each infinite fiber,
there is an Einstein-Hermitian metric $k$ of $E$.
We set $\rho = \root r \of {\det(k)/\det(h)}$ and $h' = \rho h$.
Then it is easy to see that  $\det(h') = \det(h)$ and
$\Delta(E, h) = \Delta(E, h')$. Thus by Lemma~6.2,
$$
 \align
  \Delta(E, h) - \Delta(E, k) & = \Delta(E, h') - \Delta(E, k) \\
  & = (\achern{2}{E, h'} - \achern{2}{E, k})\cdot\achern{1}{H, h_H}^{d-2} \\
  & = \frac{(d-2)!}{2\pi} \sum_{\sigma \in K_{\infty}}
  DL(E_{\sigma}, h'_{\sigma}, k_{\sigma} ; \Phi_{\sigma}).
 \endalign
$$
Here, since $DL$ is the Donaldson's Lagrangian, we get
$$
\sum_{\sigma \in S_{\infty}}
  DL(E_{\sigma}, h'_{\sigma}, k_{\sigma} ; \Phi_{\sigma}) \geq 0
$$
and the equality holds if and only if $h'$ is Einstein-Hermitian
on each infinite fiber. Thus we have our theorem.
\qed
\enddemo
\vfill\eject

\head 7. Second fundamental form
\endhead

Let $M$ be a complex manifold and $(E, h)$ a Hermitian vector bundle
on $M$. Let $0 \to S \to E \to Q \to 0$ be an exact sequence of vector
bundles. Let $h'$ and $h''$ be Hermitian metrics of $S$ and
$Q$ induced by $h$ respectively.
Let $E = S \oplus S^{\perp}$ be the orthogonal decomposition of $E$ by $h$.
Let $D(E, h)$, $D(S, h')$ and $D(Q, h'')$ be the Hermitian connections
of $(E, h)$, $(S, h')$ and $(Q, h'')$ respectively.
Moreover, let $K(E, h)$, $K(S, h')$ and $K(Q, h'')$ be the curvatures of
$(E, h)$, $(S, h')$ and $(Q, h'')$ respectively.
The Hermitian connection $D(E, h)$ has the following form:
$$
   D(E, h) = \pmatrix D(S, h') & -A^{*} \\
                      A        & D(Q, h'')
             \endpmatrix,
$$
where $A \in A^{1, 0}(\operatorname{Hom}(S, S^{\perp}))$ and
$A^{*}$ is the adjoint of $A$. $A$ is called
{\it the second fundamental form of}
$$0 \to (S, h') \to (E, h) \to (Q, h'') \to 0.$$
It is well known that
the exact sequence $0 \to S \to E \to Q \to 0$
induces the orthogonal decomposition $(E, h) = (S, h'') \oplus (Q, h'')$
if and only if $A$ vanishes identically.
We set
$$
  D_t = D(E, h) + (e^t - 1) \pmatrix 0 & 0 \\ A & 0 \endpmatrix
  \quad \hbox{and} \quad K_t = (D_t)^2.
$$
Then, by Proposition~3.28 and Lemma~4.7 of \cite{BC}, we have
$$
\multline
\tr(K(E, h)^2) - \tr(K_t^2)  = \\
2 \pi \sqrt{-1} d d^c \left(
\int_{t}^{0} \left\{
  \tr \left( K_t \cdot \pmatrix 1 & 0 \\ 0 & 0 \endpmatrix \right) +
  \tr \left( \pmatrix 1 & 0 \\ 0 & 0 \endpmatrix \cdot K_t \right) \right\} dt
\right).
\endmultline
$$
Therefore, since $\tr(K(S, h'))$ is $d$-closed,
by an easy calculation, we get
$$
\tr(K(E, h)^2) - \tr(K_t^2) =
-4 \pi \sqrt{-1} (1-e^t)d d^c \left( \tr(A^* \wedge A) \right).
$$
Thus we obtain
$$
\tr(K(E, h)^2) - \tr(K(S, h')^2 \oplus K(Q, h'')^2)
= -4 \pi \sqrt{-1} d d^c \left( \tr (A^* \wedge A) \right)
$$
because
$$
\lim_{t \to -\infty} \tr(K_t^2) = \tr(K(S, h')^2 \oplus K(Q, h'')^2).
$$
These observations show us the following lemma.

\proclaim{Lemma 7.1}
With notation being as above, we have
$$
\chern_2(E, h) - \chern_2((S, h') \oplus (Q, h''))
= d d^c \left( \frac{\sqrt{-1}}{2 \pi} \tr (A^* \wedge A) \right).
$$
In particular, by the axiomatic characterization of $\sch_2$, we get
$$
\sch_2(0 \to (S, h') \to (E, h) \to (Q, h'') \to 0)
= \frac{\sqrt{-1}}{2 \pi} \tr (A^* \wedge A)
$$
modulo $\partial(A^{0,1}) + \bar{\partial}(A^{1,0})$.
\endproclaim

\definition{Remark 7.2}
In the sense of \cite{BGS},
$$
\sch_2(0 \to (S, h') \to (E, h) \to (Q, h'') \to 0)
= -\frac{\sqrt{-1}}{2 \pi} \tr (A^* \wedge A).
$$
\enddefinition

Here we assume that $M$ is an $n$-dimensional compact K\"{a}hler manifold with
a K\"{a}hler form $\Phi$.
Let $\theta^1, \ldots, \theta^n$ be a local unitary frame of $\Omega_M^1$.
Then $\Phi = \sqrt{-1} \sum \theta^i \wedge \bar{\theta}^i$.
We set $A = \sum A_i \theta^i$. Since $A^{*} = \sum A_i^{*} \bar{\theta}^i$,
we get
$$
\align
\sqrt{-1} \tr(A^* \wedge A) \wedge \frac{\Phi^{n-1}}{(n-1)!}
& =
\sqrt{-1} \sum_{i=1, j=1}^n \tr(A_i^* \wedge A_j)
(\bar{\theta}^i \wedge \theta^j) \wedge \frac{\Phi^{n-1}}{(n-1)!} \\
& =
-\sum_{i=1}^n \tr(A_i^* \wedge A_i) \sqrt{-1}(\theta^i \wedge \bar{\theta}^i)
\wedge \frac{\Phi^{n-1}}{(n-1)!} \\
& = - | A |^2 \frac{\Phi^n}{n!},
\endalign
$$
which implies that
$$
\int_M \sch_2(0 \to (S, h') \to (E, h) \to (Q, h'') \to 0) \wedge
\frac{\Phi^{n-1}}{(n-1)!} = -\frac{1}{2 \pi} ||A||^2.
$$
Thus we have the following proposition.

\proclaim{Proposition 7.3}
Let $K$ be an algebraic number field and $O_K$ the ring of integers.
We denote by $K_{\infty}$ the set of
all embeddings of $K$ into $\CC$.
Let $f : X \longrightarrow \Spec(O_K)$ be a regular arithmetic variety with
$\dim X = d \geq 2$, and $(H, h_H)$ a Hermitian line bundle on $X$
such that, for each $\sigma \in K_{\infty}$,
$c_1(H_{\sigma}, h_{H_{\sigma}})$ gives a K\"{a}hler form $\Phi_{\sigma}$ on
an infinite fiber $X_{\sigma}$.
Let $0 \to S \to E \to Q \to 0$ be an
exact sequence of torsion free sheaves such that
each torsion free sheaf is locally free on the generic fiber.
Let $h$ be a Hermitian metric of $E$, and $h'$ and $h''$ Hermitian metrics
of $S$ and $Q$ induced by $h$ respectively.
Then, we have the following:
$$
\left( \ach_2(E, h) - \ach_2((S, h') \oplus (Q, h'')) \right)
\cdot \achern{1}{H, h_H}^{d-2}
= - \frac{(d-2)!}{2 \pi} \sum_{\sigma \in K_{\infty}} ||A_{\sigma}||^2,
$$
$$
\left( \achern{2}{E, h} - \achern{2}{(S, h') \oplus (Q, h'')} \right)
\cdot \achern{1}{H, h_H}^{d-2}
= \frac{(d-2)!}{2 \pi} \sum_{\sigma \in K_{\infty}} ||A_{\sigma}||^2.
$$
\endproclaim

\proclaim{Corollary 7.4}
Let $K$ be an algebraic number field and $O_K$ the ring of integers.
We denote by $K_{\infty}$ the set of
all embeddings of $K$ into $\CC$.
Let $f : X \longrightarrow \Spec(O_K)$ be a regular arithmetic variety with
$\dim X = d \geq 2$, and $(H, h_H)$ a Hermitian line bundle on $X$
such that, for each $\sigma \in K_{\infty}$,
$c_1(H_{\sigma}, h_{H_{\sigma}})$ gives a K\"{a}hler form $\Phi_{\sigma}$ on
an infinite fiber $X_{\sigma}$. Let
$$ 0 = E_0 \subset E_1 \subset \cdots \subset E_{l-1} \subset E_l
$$
be a filtration of torsion free sheaves on $X$ such that
\roster
\item "(i)" $E_i$ is locally free on the generic fiber
            for every $1 \leq i \leq l$, and that

\item "(ii)" $E_i/E_{i-1}$ is torsion free and locally free on the
             generic fiber for every $1 \leq i \leq l$.
\endroster
Let $h_l$ be a Hermitian metric of $E_l$ and $h_i$ the induced metric of
$E_i$ by $h_l$. Let $Q_i = E_i/E_{i-1}$ and $k_i$ the quotient metric of
$Q_i$ induced by $h_i$. Then, we have
$$
\bigl\{ \achern{2}{E_l, h_l} -
 \achern{2}{(Q_1, k_1) \oplus \cdots \oplus (Q_l, k_l)} \bigr\}
\cdot \achern{1}{H, h_H}^{d-2} \geq 0.
$$
Further, the equality of the above inequality holds if and only if
$(E_l, h_l)$ is isometric to $(Q_1, h_1) \oplus \cdots \oplus (Q_l, h_l)$
on each infinite fiber.
\endproclaim

\demo{Proof}
We will prove the following inequality inductively.
$$
\bigl\{ \achern{2}{E_l, h_l} -
 \achern{2}{(E_i, h_i) \oplus (Q_{i+1}, k_{i+1})
 \oplus \cdots \oplus (Q_l, k_l)} \bigr\}
\cdot \achern{1}{H, h_H}^{d-2} \geq 0.
$$
In the case where $i = l-1$, the above inequality is an immediate consequence
of Proposition~7.3.
Here we consider an exact sequence:
$$ 0 \to E_{i-1} \to E_i \to Q_i \to 0. $$
By Proposition~7.3, we have
$$
\bigl\{ \achern{2}{E_i, h_i} - \achern{2}{(E_{i-1}, h_{i-1})
\oplus (Q_i, k_i)} \bigr\} \cdot \achern{1}{H, h_H}^{d-2} \geq 0.
$$
Therefore, we obtain
$$
 \multline
 \bigl\{
 \achern{2}{(E_i, h_i) \oplus (Q_{i+1}, k_{i+1})
 \oplus \cdots \oplus (Q_l, k_l)} - \\ \quad
 \achern{2}{(E_{i-1}, h_{i-1}) \oplus (Q_i, k_i) \oplus (Q_{i+1}, k_{i+1})
 \oplus \cdots \oplus (Q_l, k_l)} \bigr\}
 \cdot \achern{1}{H, h_H}^{d-2} \geq 0.
 \endmultline
$$
Hence, combining hypothesis of induction, we get
$$
\bigl\{ \achern{2}{E_l, h_l} -
 \achern{2}{(E_{i-1}, h_{i-1}) \oplus (Q_i, k_i) \oplus
 \cdots \oplus (Q_l, k_l)} \bigr\}
\cdot \achern{1}{H, h_H}^{d-2} \geq 0.
$$
Thus we have our corollary.
\qed
\enddemo
\vfill\eject

\head 8. Proof of the main theorem
\endhead

First of all, we will prepare the following two lemmas.

\proclaim{Lemma 8.1}
Let $k$ an algebraically closed field and $K$ an extension field of $k$.
Let $X$ be a normal projective variety over $k$, $H$ an ample line bundle
on $X$, and $E$ a vector bundle on $X$.
Then $E$ is $H$-stable (resp. $H$-semistable) if and only if
$E \otimes_k \overline{K}$ is $H \otimes_k \overline{K}$-stable
(resp. $H \otimes_k \overline{K}$-semistable).
\endproclaim

\demo{Proof}
Clearly, if
$E \otimes_k \overline{K}$ is $H \otimes_k \overline{K}$-stable
(resp. $H \otimes_k \overline{K}$-semistable),
then $E$ is $H$-stable (resp. $H$-semistable).

We assume that $E \otimes_k \overline{K}$ is not
$H \otimes_k \overline{K}$-stable
(resp. not $H \otimes_k \overline{K}$-semistable).
Then there is a subsheaf $F$ of $E \otimes_k \overline{K}$
with $\mu(F) \geq \mu(E)$ (resp. $\mu(F) > \mu(E)$).
Here we can take elements $x_1, \cdots, x_n$ of $\overline{K}$ such that
$F$ is defined over $k(x_1, \ldots, x_n)$.
We consider a normal projective variety $Y$ over $k$ such that
$k(Y) = k(x_1, \ldots, x_n)$.
We set $Z = X \times_k Y$ and let $p : Z \to X$ be the projection.
By the assumption, for the generic point $\eta \in Y$,
$\rest{p^*(E)}{X \times \{ \bar{\eta} \}}$ is not $H$-stable
(resp. not $H$-semistable). Therefore, by \cite{Ma1},
$\rest{p^*(E)}{X \times \{ y \}}$ is not $H$-stable
(resp. not $H$-semistable) for all closed point $y$ of $Y$.
Thus $E$ is not $H$-stable (resp. not $H$-semistable).
\qed
\enddemo

\proclaim{Lemma 8.2}
Let $k$ be a Galois extension over $\QQ$ with a Galois group
$G = \Gal(k/\QQ)$.
Let $X$ be a smooth projective variety over $k$, $H$ an ample line bundle
on $X$, and $E$ a vector bundle on $X$.
Then if $E \otimes_k \overline{\QQ}$ is stable (resp. semistable),
then, for every $\sigma \in G$, $E \otimes_k^{\sigma} \overline{\QQ}$ is
also stable (resp. semistable), where $E \otimes_k^{\sigma} \overline{\QQ}$
is a tensor product with using an embedding $\sigma : k \to \overline{\QQ}$.
\endproclaim

\demo{Proof}
Assume that $E \otimes_k^{\sigma} \overline{\QQ}$ is not stable
(resp. not semistable). Then there is a subsheaf $F$
of $E \otimes_k^{\sigma} \overline{\QQ}$ such that
$\mu(F) \geq \mu(E)$ (resp. $\mu(F) > \mu(E)$).
We may assume that $F$ is defined over a field $k'$ such that
$k \subset k'$ and $k'$ is a Galois extension over $\QQ$.
We take an element $\sigma'$ of $\Gal(k'/\QQ)$ with
$\rest{\sigma'}{k} = \sigma$.
Here we give a right $k'$-module structure to $(E \otimes_k^{\sigma} k')$
by $(e \otimes a) \cdot \lambda = e \otimes a \sigma'(\lambda)$.
Then if we consider a correspondence $e \otimes a \otimes b
\rightsquigarrow e \otimes {\sigma'}^{-1}(a) \sigma'(b)$,
we can easily to see that
$(E \otimes_k^{\sigma} k') \otimes_{k'}^{{\sigma'}^{-1}} k' \simeq
E \otimes_k k'$. Moreover, if we give a right $k'$-module structure
to $(E \otimes_k^{\sigma} k') \otimes_{k'}^{{\sigma'}^{-1}} k'$
by $(e \otimes a \otimes b) \cdot \lambda =
e \otimes a \otimes b {\sigma'}^{-1}(\lambda)$, then the above is
an isomorphism as $k'$-modules.
Thus $F \otimes_{k'}^{{\sigma'}^{-1}} k'$
is a subsheaf of $E \otimes_k k'$.
On the other hand, since the intersection number does not change by
an extension of the grand field,
we have $\mu(F) = \mu(F \otimes_{k'}^{{\sigma'}^{-1}} k')$.
This is a contradiction.
\qed
\enddemo

Let us start the proof of the main theorem.
We follow steps in Introduction.

\subhead 8.3 Step 1 \endsubhead
We try to reduce our theorem to the case where
$E$ is poly-stable on each infinite fiber.
Since $E \otimes_K \overline{\QQ}$ is semistable,
there is a Jordan-H\"{o}lder filtration of $E \otimes_K \overline{\QQ}$:
$$
0 = E_0 \subset E_1 \subset \cdots \subset E_{l-1} \subset E_l =
E \otimes_K \overline{\QQ}
$$
such that $E_{i}/E_{i-1}$ is stable for all $1 \leq i \leq l$
and $\mu(E_1/E_0) = \mu(E_2/E_1) = \cdots = \mu(E_{l}/E_{l-1})$.
Considering a base change of $O_K$, we may assume that
$E_i$ is defined over $K$ and $K$ is a Galois extension over $\QQ$.
Moreover, if we take a suitable birational change of $X$, we may assume that
$X$ is regular, $E_i$ is defined over $X$ and that
$E_i/E_{i-1}$ is torsion free. We set $Q_i = E_i/E_{i-1}$.
We give a Hermitian metric $h_i$ to each $E_i$ induced by $h$.
Here we consider an exact sequence:
$$
   0 \to E_{i-1} \to E_i \to Q_i \to 0.
$$
Let $k_i$ be the quotient metric of $Q_i$ by the above exact sequence, and
$Q_i^{\vee\vee}$ the double dual of $Q_i$.
Lemma~8.1 and Lemma~8.2 imply
$Q_1^{\vee\vee} \oplus \cdots \oplus Q_l^{\vee\vee}$ is poly-stable
on each infinite fiber. Thus by hypothesis of reduction,
we have
$$
 \achern{2}{(Q_1^{\vee\vee}, k_1) \oplus \cdots \oplus (Q_l^{\vee\vee}, k_l)}
 - \frac{r-1}{2r}
 \achern{1}{(Q_1^{\vee\vee}, k_1) \oplus \cdots \oplus (Q_l^{\vee\vee}, k_l)}^2
 \geq 0.
$$
Corollary~7.4 implies that
$$
 \achern{2}{E, h} \geq
 \achern{2}{(Q_1, k_1) \oplus \cdots \oplus (Q_l, k_l)}.
$$
On the other hand, clearly we have
$$
\achern{2}{(Q_1, k_1) \oplus \cdots \oplus (Q_l, k_l)} \geq
\achern{2}{(Q_1^{\vee\vee}, k_1) \oplus \cdots \oplus (Q_l^{\vee\vee}, k_l)}
$$
and
$$
 \achern{1}{E, h} =
 \achern{1}{(Q_1^{\vee\vee}, k_1) \oplus \cdots \oplus (Q_l^{\vee\vee}, k_l)}.
$$
Thus we get Step 1.

\subhead 8.4 Step 2 \endsubhead
By Theorem~6.3,
we may assume that $h$ is Einstein-Hermitian on each infinite fiber.

\subhead 8.5 Step 3 \endsubhead
Let $\pi : Y = \PP(E) \longrightarrow X$ be the projective bundle of $E$ and
$\OO_Y(1)$ the tautological line bundle.
Let $\Cal{E} : 0 \to F \to \pi^*E \to \OO_Y(1) \to 0$ be the canonical
exact sequence.
We give the canonical Hermitian metrics on
$F$, $\pi^*E$ and $\OO_Y(1)$ induced by the Hermitian metric of $E$.
We set $L = \OO_Y(r) \otimes \pi^*(\det E)^{-1}$.
Here we consider
$$
\deg \left(\achern{1}{\OO_Y(1)} - \frac{1}{r}\pi^*(\achern{1}{E})
\right)^{r+1}.
$$
We set
$$
\Phi =
\sum_{\sigma \in S_{\infty}}\sum_{i=1}^r (-1)^i
\tilde{c}_i(\Cal{E}_{\sigma}) c_1(\OO_{Y_{\sigma}}(1))^{r-i}.
$$
Then by (1.9.1)
$$
   \achern{1}{\OO_{Y}(1)}^r - \pi^*(\achern{1}{E}) \achern{1}{\OO_{Y}(1)}^{r-1}
   + \pi^*(\achern{2}{E}) \achern{1}{\OO_{Y}(1)}^{r-2}
   = a(\Phi).
$$
Thus we have
$$
\align
  \left(\achern{1}{\OO_Y(1)} - \frac{1}{r}\pi^*(\achern{1}{E})\right)^{r+1} = &
  \frac{r-1}{2r} \pi^*(\achern{1}{E}^2)\achern{1}{\OO_Y(1)}^{r-1} \\
  & \quad - \pi^*(\achern{2}{E})\achern{1}{\OO_Y(1)}^{r-1} \\
  & \quad
  + \left( \achern{1}{\OO_Y(1)}-\frac{1}{r}\pi^*(\achern{1}{E}) \right)
a(\Phi).
\endalign
$$
Since
$$
  \pi_*(\pi^*(\achern{1}{E}^2)\achern{1}{\OO_Y(1)}^{r-1}) = \achern{1}{E}^2
$$
and
$$
\pi_*(\pi^*(\achern{2}{E})\achern{1}{\OO_Y(1)}^{r-1}) = \achern{2}{E}
$$
by projection formula, we have
$$
\align
  \pi_* \left(\achern{1}{\OO_Y(1)} -
             \frac{1}{r}\pi^*(\achern{1}{E})\right)^{r+1} & =
  \frac{r-1}{2r} (\achern{1}{E}^2) -
                 (\achern{2}{E}) + \\
  & \quad \pi_* \left(
      \left( \achern{1}{\OO_Y(1)}-\frac{1}{r}\pi^*(\achern{1}{E}) \right)
         a(\Phi) \right).
\endalign
$$
On the other hand, since
$$
 \deg \left(
 \left( \achern{1}{\OO_Y(1)}-\frac{1}{r}\pi^*(\achern{1}{E}) \right) a(\Phi)
      \right)
$$
is equal to
$$
   \sum_{\sigma \in S_{\infty}} \int_{Y_{\sigma}} \left\{
   \left( c_1(\OO_{Y_{\sigma}}(1))-\frac{1}{r}\pi^*(c_1(E_{\sigma})) \right)
   \sum_{i=1}^r (-1)^i
   \tilde{c}_i(\Cal{E}_{\sigma}) c_1(\OO_{Y_{\sigma}}(1))^{r-i} \right\},
$$
Lemma~5.1 implies
$$
 \deg \left(
 \left( \achern{1}{\OO_Y(1)}-\frac{1}{r}\pi^*(\achern{1}{E}) \right) a(\Phi)
 \right) = 0.
$$
Therefore, we get
$$
\deg \left(\achern{1}{\OO_Y(1)} - \frac{1}{r}\pi^*(\achern{1}{E}) \right)^{r+1}
=  \frac{r-1}{2r} (\achern{1}{E}^2) - (\achern{2}{E}).
$$
On the the other hand,
$$
   (L^{r+1}) = r^{r+1}
   \deg \left(
      \achern{1}{\OO_Y(1)} - \frac{1}{r}\pi^*(\achern{1}{E})
        \right)^{r+1}.
$$
Thus we have
$$
\hbox{
$(L^{r+1}) \leq 0$ $\Longrightarrow$
${\displaystyle \frac{r-1}{2r} \achern{1}{E,h}^2 \leq \achern{2}{E,h}}$.}
$$
Hence it is sufficient to show that $(L^{r+1}) \leq 0$.

\subhead 8.6 Step 4 \endsubhead
Let $(N, h)$ be a Hermitian line bundle on $X$ such that
$N$ is ample, $\deg(N_K) > 2 g(X_K) - 2$ and that $L \otimes \pi^* N$ is ample.
By the arithmetic Grothendieck-Riemann-Roch theorem (1.11.1),
we have
$$
   \chi_{L^2}(L^n \otimes \pi^*N) +
   \frac{1}{2} \sum_{\sigma \in S_{\infty}}
               \tau((L^n \otimes \pi^*N)_{\sigma})
   = \deg(\ach(L^n \otimes \pi^*N) \cdot \atodd(T_{Y/S})).
$$
Clearly
$$
   \lim_{n\to\infty}
   \frac{\deg(\ach(L^n \otimes \pi^*N) \cdot \atodd(T_{Y/S}))}{n^{r+1}}
   = \frac{1}{(r+1)!} (L^{r+1}).
$$
Thus it is sufficient to show that
\roster
\item "(a)" ${\displaystyle
\sum_{\sigma \in S_{\infty}}
\tau((L^n \otimes \pi^*N)_{\sigma}) \leq O(n^r \log n)}$.

\item "(b)" $\chi_{L^2}(L^n \otimes \pi^*N) \leq O(n^r \log n)$.
\endroster

\subhead 8.7 Step 5 \endsubhead
If we denote by $H_{\sigma}$ the Hermitian form of $L_{\sigma}$
for all $\sigma \in S_{\infty}$,
$H_{\sigma}(y)$ is positive semi-definite and $\rank H_{\sigma}(y) \geq r-1$
for all $y \in Y_{\sigma}$ as in p.89-p.90 of \cite{Ko}.
Thus by Corollary~2.4 we have (a).

\subhead 8.8 Step 6 \endsubhead
Since $\deg(N) > 2 g(X) -2$ and $\Sym^{rn}(E)$ is semistable
on the generic fiber, $H^1(Y, L^n \otimes \pi^*N)$ is torsion module.
Thus it is sufficient to see that
$$
   \deg_{L^2}(H^0(L^n \otimes \pi^*N)) \leq O(n^r \log n).
$$
Since $L$ is nef and $(L^r) = 0$ on the generic fiber,
we have
$$
 \rank H^0(Y, L^{n} \otimes \pi^*N) \leq O(n^{r-1}).
$$
Hence Theorem~4.1 implies that
$$
   \deg_{L^2}(H^0(Y, L^{n} \otimes \pi^*N)) \leq O(n^{r} \log n).
$$
Thus we get (b), which completes the proof of the main theorem.
\qed

\proclaim{Corollary 8.9}
Let $K$ be an algebraic number field, $O_K$ the ring of integers of $K$,
and $f : X \longrightarrow \Spec(O_K)$ an arithmetic surface.
Let $E$ be a torsion free sheaf and $h$ a Hermitian metric on $X$.
If $X$ is regular and $E_{\overline{\QQ}}$ is semistable on
the geometric generic fiber $X_{\overline{\QQ}}$ of $f$,
then we have an inequality
$$
\achern{2}{E, h} - \frac{r-1}{2r} \achern{1}{E, h}^2 \geq 0,
$$
where $r = \rank E$.
\endproclaim

\demo{Proof}
Let $E^{\vee\vee}$ be the double dual of $E$. Since $X$ is a regular scheme
of dimension $2$, $E^{\vee\vee}$ is locally free and $T = E^{\vee\vee}/E$
is a finite module. Moreover, $\achern{1}{E, h} = \achern{1}{E^{\vee\vee}, h}$
and $\achern{2}{E, h} = \achern{2}{E^{\vee\vee}, h} + \log \#(T)$.
Thus we have our corollary.
\qed
\enddemo
\vfill\eject

\head 9. Arithmetic second Chern character of semistable vector bundles
\endhead

Let $K$ be an algebraic number field and $O_K$ the ring of integers.
Let $f : X \longrightarrow \Spec(O_K)$ be a regular arithmetic surface and
$E$ a torsion free sheaf on $X$.
We denote by $\operatorname{Herm}(E)$ the set of all Hermitian metrics of $E$.
We set
$$
  \achsup{2}{E} = \sup_{h \in \operatorname{Herm}(E)} \ach_2(E, h)
  \quad (\in (-\infty, \infty]).
$$
First of all, we will see several properties of $\achsup{2}{E}$ when
$E$ is semistable on the geometric generic fiber.

\proclaim{Proposition 9.1}
With being notation as above, we assume that
$E_{\overline{\QQ}}$ is semistable and $\deg(E_{\overline{\QQ}}) = 0$.
Then, we have the following.
\roster
 \item "(1)"
       $\ach_2(E, h) \leq 0$ for every $h \in \operatorname{Herm}(E)$.
       Moreover, the equality holds if and only if
       $\achern{1}{E, h}^2 = \achern{2}{E, h} = 0$.

 \item "(2)"
       $\achsup{2}{E} \leq 0$.

 \item "(3)"
       For $h \in \operatorname{Herm}(E)$, $\achsup{2}{E} = \ach_2(E, h)$
       if and only if $h$ is Einstein-Hermitian.

 \item "(4)"
       If $\rank E = 1$, then
       $\achsup{2}{E} \leq -[K : \QQ] \operatorname{Height}([\rest{E}{K}])$,
       where $\operatorname{Height}$ is the N\'{e}ron-Tate height function
       on $\Pic^0(Y)$. Moreover, the equality holds if and only if
       $E$ is locally free and
       $\deg(\rest{E}{C}) = 0$ for all vertical curves $C$ on $X$.
 \item "(5)"
       If there is a filtration of $E$:
       $0 = E_0 \subset E_1 \subset \cdots \subset E_{l-1} \subset E_l = E$
       such that $E_i/E_{i-1}$ is torsion free and
       $\deg(\rest{(E_i/E_{i-1})}{K}) = 0$ for every $i$,
       then we have $$ \achsup{2}{E} = \sum_{i=1}^l \achsup{2}{E_i/E_{i-1}}. $$

 \item "(6)"
       We assume that $E$ is locally free.
       Let $K'$ be a finite extension field of $K$, $X'$
       a desingularization of $X \times_{O_K} O_{K'}$, and $g : X' @>>> X$
       the induced morphism. Then, we have
       $$ \achsup{2}{g^*(E)} = [K' : K] \achsup{2}{E}. $$
\endroster
\endproclaim

\demo{Proof}
(1) First, we will see that $\achern{1}{E, h}^2 \leq 0$.
We set $$(E', h') = (\det(E), \det(h)) \oplus (\OO_X, \operatorname{can}_K).$$
Then, since
$E'$ is semistable, by Corollary~8.9, we get
$$
      \achern{1}{E, h}^2 = \achern{1}{E', h'}^2 \leq
      4 \achern{2}{E', h'} = 0.
$$
Hence, by virtue of Corollary~8.9,
$$
\align
\ach_2(E, h) & = \frac{1}{2} \achern{1}{E, h}^2 - \achern{2}{E, h} \\
& \leq \frac{1}{2r} \achern{1}{E, h}^2 - \left(
\achern{2}{E, h} - \frac{r-1}{2r} \achern{1}{E, h}^2 \right) \leq 0.
\endalign
$$

\medskip
(2) This is a direct consequence of (1).

\medskip
(3) We fix a Hermitian metric $k$ of $E$.
Then, since $\deg(L_K) = 0$, by Lemma~6.1,
$$
\ach_2(E, h) - \ach_2(E, k) = -\frac{1}{2\pi} \sum_{\sigma \in K_{\infty}}
 DL(E_{\sigma}, h_{\sigma}, k_{\sigma} ; \Phi_{\sigma}).
$$
Thus, $\ach_2(E, h) = \achsup{2}{E}$
if and only if $DL(E_{\sigma}, {h}_{\sigma}, k_{\sigma} ; \Phi_{\sigma})$
gives the absolute minimal value for every $\sigma \in K_{\infty}$.
Thus, we obtain (3).

\medskip
(4) Let $h$ be an Einstein-Hermitian metric of $E$. By (3), we have
$\achsup{2}{E} = \ach_2(E, h)$.
Let $E^{\vee\vee}$ be the double dual of $E$. Then,
$\ach_2(E, h) \leq \ach_2(E^{\vee\vee})$ and the equality holds if and
only if $E = E^{\vee\vee}$. So we may assume that $E$ is locally free.
Here, we need the following lemma.

\proclaim{Lemma 9.2}
Let $(L,  h)$ a Hermitian line bundle on $X$.
If $h$ is Einstein-Hermitian and $\deg(\rest{L}{C}) = 0$ for
all vertical curves $C$ on $X$, then
$$
 \ach_2(L, h) = - [K : \QQ] \operatorname{Height}([L_K]).
$$
\endproclaim

\demo{Proof} For example, see \cite{Fa1} and \cite{Hr}. \qed
\enddemo

\medskip
Let $\{ q_1, \ldots, q_l \}$ be the set of all critical values of $f$.
Let $\Supp(X_{q_t}) = C^{(t)}_1 + \cdots + C^{(t)}_{e_t}$ be the irreducible
decomposition of the fiber $X_{q_t}$.
Then, since $\deg(\rest{E}{K}) = 0$,
we can find rational numbers $\{ a^{(t)}_i \}$ such that
$$
    (E \cdot C^{(t)}_j) +
    (\sum_{i=1}^{e_t} a^{(t)}_i C^{(t)}_i \cdot C^{(t)}_j) = 0
$$
for all $t$ and $j$.
Let $n$ be a positive integer such that $n \cdot a^{(t)}_i \in \ZZ$
for all $t$ and $i$.
We set
$$
  L = E^n \otimes \OO_X(\sum_{t=1}^l \sum_{i=1}^{e_t} n a^{(t)}_i C^{(t)}_i).
$$
Then, due to Lemma~9.2,
we have
$$
   \ach_2(L, h^n) = -[K : \QQ] \operatorname{Height}([L_K])
                  = - n^2 [K : \QQ] \operatorname{Height}([E_K]).
$$
On the other hand, since $\deg(\rest{L}{C}) = 0$ for all
vertical curves $C$ on $X$, we have
$$
 \ach_2(L, h^n) + \frac{n^2}{2}
 (\sum_{t=1}^l \sum_{i=1}^{e_t} a^{(t)}_i C^{(t)}_i)^2 =
 n^2 \ach_2(E, h).
$$
Let $\sum_{t=1}^l \sum_{i=1}^{e_t} a^{(t)}_i C^{(t)}_i = D_{+} - D_{-}$
be the decomposition of $\sum_{t=1}^l \sum_{i=1}^{e_t} a^{(t)}_i C^{(t)}_i$
such that $D_{+}$ and $D_{-}$ are effective and that
$\Supp(D_{+})$ and $\Supp(D_{-})$ have no common component.
Then,
$$
(\sum_{t=1}^l \sum_{i=1}^{e_t} a^{(t)}_i C^{(t)}_i)^2
= (D_{+}^2) + (D_{-}^2) - 2(D_{+} \cdot D_{-}) \leq 0.
$$
Moreover, the equality holds if and only if
$$
 (D_{+} \cdot C) = (D_{-} \cdot C) = 0
$$
for all vertical curves $C$ on $X$. Thus we have (4).

\medskip
(5) Let $h$ be a Hermitian metric of $E$.
Let $h_i$ be the sub-metric of $E_i$ induced by $h$
and $k_i$ the quotient metric of $E_i/E_{i-1}$ induced by $h_i$.
Then, by Corollary~7.4, we have
$$
\align
\ach_2(E, h) & \leq
\ach_2(E_1/E_0, k_1) + \cdots + \ach_2(E_{l}/E_{l-1}, k_l) \\
& \leq
\achsup{2}{E_1/E_0} + \cdots + \achsup{2}{E_{l}/E_{l-1}}
\endalign
$$
Therefore, we have
$$
\achsup{2}{E} \leq \sum_{i=1}^l \achsup{2}{E_i/E_{i-1}}.
$$

In order to consider the converse inequality of the above,
we need the following lemma.

\proclaim{Lemma 9.3}
Let $0 \to S \to E \to Q \to 0$ an exact sequence of torsion free sheaves
on $X$. Let $h'$ and $h''$ be Hermitian metrics of $S$ and $Q$ respectively.
If $\deg(S_K) = \deg(E_K) = \deg(Q_K) = 0$, then there is a
family $\{ h_t \}_{t \in \RR}$ of Hermitian metrics of $E$ with
$$
  \lim_{t \to \infty} \ach_2(E, h_t) = \ach_2((S, h') \oplus (Q, h'')).
$$
\endproclaim

\demo{Proof}
For $\sigma \in K_{\infty}$, let $E_{\sigma} = S_{\sigma} \oplus P_{\sigma}$
be a decomposition of $E_{\sigma}$ as $C^{\infty}$-vector bundles.
Then, since $P_{\sigma}$ is isomorphic to $Q_{\sigma}$,
using the above decomposition, we define a Hermitian metric $h_t$ of $E$
by $e^t \cdot h' \oplus h''$.
Let $A_t$ be the second fundamental form of
$$ 0 \to (S, e^t \cdot h') \to (E, h_t) \to (Q, h'') \to 0.$$
If we denote $A_0$ by $A$, then it is easy to see that
$A_t = A$ and $(A_t)^* = e^{-t} A^*$.
Thus, by Proposition~7.3, we have
$$
\ach_2(E, h_t) - \ach_2((S, e^t \cdot h') \oplus (Q, h''))
= - \frac{e^{-t}}{2 \pi} \sum_{\sigma \in K_{\infty}} ||A_{\sigma}||^2.
$$
On the other hand, since $\deg(S_K) = 0$, by an easy calculation,
we obtain $\ach_2(S, e^t \cdot h') = \ach_2(S, h')$.
Thus, we get our lemma.
\qed
\enddemo

\medskip
Let us start to prove the inequality
$$
\achsup{2}{E} \geq \sum_{i=1}^l \achsup{2}{E_i/E_{i-1}}.
$$
by induction on $l$.
Let $h'$ and $h''$ be arbitrary Hermitian metrics of $E_{l-1}$ and $Q_l$
respectively. Then, by Lemma~9.3,
there is a family $\{ h_t \}_{t \in \RR}$ of
Hermitian metric of $E$ such that
$$
\lim_{t \to \infty} \ach_2(E, h_t) = \ach_2(E_{l-1}, h') + \ach_2(Q_l, h'').
$$
Therefore, we have
$$
\achsup{2}{E} \geq \ach_2(E_{l-1}, h') + \ach_2(Q_l, h'')
$$
Since $h'$ and $h''$ are arbitrary, it follows that
$$
\achsup{2}{E} \geq \achsup{2}{E_{l-1}} + \achsup{2}{Q_l}
$$
Thus, by hypothesis of induction, we have (5).

\medskip
(6) Let $k$ be a fixed Hermitian metric of $E$.
For each $\sigma \in K_{\infty}$, we set
$$
  \overline{DL}(E_{\sigma}) = \inf_{h_{\sigma}}
  DL(E_{\sigma}, h_{\sigma}, k_{\sigma} ; \Phi_{\sigma}),
$$
where $h_{\sigma}$ runs over all Hermitian metrics of $E_{\sigma}$.
Let $\{ h_n \}$ be a sequence of Hermitian metrics of $E$ with
$$
       \lim_{n \to \infty} \ach_2(E, h_n) = \achsup{2}{E}.
$$
Since
$$
\ach_2(E, h_n) - \ach_2(E, k) = -\frac{1}{2\pi} \sum_{\sigma \in K_{\infty}}
 DL(E_{\sigma}, (h_n)_{\sigma}, k_{\sigma} ; \Phi_{\sigma}),
$$
for each $\sigma \in K_{\infty}$, we have
$$
  \overline{DL}(E_{\sigma}) = \lim_{n \to \infty}
  DL(E_{\sigma}, (h_n)_{\sigma}, k_{\sigma} ; \Phi_{\sigma}).
$$
Therefore, we obtain
$$
  \overline{DL}(g^*(E)_{\sigma'}) = \lim_{n \to \infty}
  DL(E_{\rest{\sigma'}{K}}, (h_n)_{\rest{\sigma'}{K}},
     k_{\rest{\sigma'}{K}} ; \Phi_{\rest{\sigma'}{K}})
$$
for all $\sigma' \in K'_{\infty}$ because
$\overline{DL}(g^*(E)_{\sigma'}) = \overline{DL}(E_{\rest{\sigma'}{K}})$.
It follows that
$$
\align
\achsup{2}{g^*(E)} & = \lim_{n \to \infty} \ach_2(g^*(E, h_n)) \\
                   & = \lim_{n \to \infty} [K' : K] \ach_2(E, h_n) \\
                   & = [K' : K] \achsup{2}{E}
\endalign
$$
Thus we have (6).
\qed
\enddemo

\proclaim{Corollary 9.4}
Let $f : X \longrightarrow \Spec(O_K)$ be a regular arithmetic surface and
$E$ a torsion free sheaf on $X$ such $E_{\overline{\QQ}}$ is semistable
and $\deg(E_{K}) = 0$. If there is a filtration of $E$:
$0 = E_0 \subset E_1 \subset \cdots \subset E_{r-1} \subset E_r = E$
such that $E_i/E_{i-1}$ is torsion free of rank $1$ and
$\deg(\rest{(E_i/E_{i-1})}{K}) = 0$ for every $i$,
then we have
$$
   \achsup{2}{E} \leq -[K : \QQ] \sum_{i=1}^r
   \operatorname{Height}((E_i/E_{i-1})_K).
$$
Moreover, the equality holds if and only if
$E_i/E_{i-1}$ is locally free
and $\deg(\rest{(E_i/E_{i-1})}{C}) = 0$ for all vertical curves $C$ on $X$.
\endproclaim

\demo{Proof}
By (5) of Proposition~9.1, we have
$$ \achsup{2}{E} = \sum_{i=1}^r \achsup{2}{E_i/E_{i-1}}. $$
Moreover, due to (4) of Proposition~9.1,
for every $i$,
$$
        \achsup{2}{E_i/E_{i-1}} \leq
        -[K : \QQ] \operatorname{Height}((E_i/E_{i-1})_K),
$$
and the equality holds if and only if $E_i/E_{i-1}$ is locally free
and $\deg(\rest{(E_i/E_{i-1})}{C}) = 0$ for all vertical curves $C$ on $X$.
Thus we have our corollary.
\qed
\enddemo

\definition{Definition 9.5}
Let $Y$ be a smooth algebraic curves over an algebraic number field $K$.
Let $\SStable_{Y/K}(r, d)$ be the
moduli scheme of semistable vector bundles on $Y$
with rank $r$ and degree $d$ (cf. \cite{Ma2}).
\enddefinition

We would like to consider $\SStable_{Y/K}(r, 0)$
when $Y$ is an elliptic curve. For this purpose, we need the following lemma.

\proclaim{Lemma 9.6}
Let $C$ be a smooth projective curve of genus 1 over an algebraically
closed field and $E$ a rank $r$ semistable vector bundle on $X$ with
$\deg(E) = 0$.
Then, there is a filtration of $E$:
$$
  0 = E_0 \subset E_1 \subset \cdots \subset E_{r-1} \subset E_r = E
$$
such that $E_i/E_{i-1}$ is a locally free of rank $1$ and
$\deg(E_i/E_{i-1}) = 0$.
\endproclaim

\demo{Proof}
We prove this proposition by induction on $\rank E$.
Let $Q$ be a minimal quotient line bundle of $E$ and
$S$ the kernel of $E \to Q$. Then, by \cite{MS},
$$
\deg(Q) - \frac{\deg(S)}{r-1} = \deg(Q) \cdot \frac{r}{r-1} \leq 1.
$$
Since $\deg(Q)$ is non-negative integer, we have $\deg(Q) = 0$.
Thus $S$ is semistable and $\deg(S) = 0$. Therefore,
$S$ has a desired filtration by hypothesis of induction.
Hence, we obtain our lemma.
\qed
\enddemo

\bigskip
The above lemma shows us that
$$
\SStable_{Y/K}(r, 0)(\overline{\QQ}) =
\overbrace{\Pic^0(Y)({\overline{\QQ}}) \times \cdots \times
\Pic^0(Y)({\overline{\QQ}})}^{\hbox{$r$ times}}/\frak{S}_r,
$$
where $\frak{S}_r$ is the $r$th symmetric group.
Thus, $\SStable_{Y/K}(r, 0)(\overline{\QQ})$ has the canonical height
function by using N\'{e}ron-Tate height, that is,
$$
  \operatorname{Height}(e) := \operatorname{Height}(l_1) + \cdots +
  \operatorname{Height}(l_r)
$$
for an element $e = (l_1, \cdots, l_r)$ of
$\SStable_{Y/K}(r, 0)(\overline{\QQ})$.
In terms of this height function $\operatorname{Height}$
of $\SStable_{Y/K}(r, 0)$, we have

\proclaim{Corollary~9.7}
Let $f : X \to \Spec(O_K)$ be a regular arithmetic surface with
the genus of the generic fiber being one, and $E$ a torsion free sheaf
of rank $r$ on $X$. Then, we have
$$
     \achsup{2}{E} \leq -[K : \QQ] \operatorname{Height}(\left[E_K\right]),
$$
where $\left[E_K\right]$ is the class of $E_K$ in
$\SStable_{X_K/K}(r, 0)$.
Moreover, the equality holds if and only if there is a finite field extension
$K'$ of $K$ with the following properties:
\roster
 \item "(a)"
       Let $X'$ be a desingularization of $X \times_{O_K} O_{K'}$
       and $g : X' \to X$ the induced morphism.

 \item "(b)"
       There is a filtration of $g^*(E)$ :
       $0 = E'_0 \subset E'_1 \subset \cdots
                 \subset E'_{r-1} \subset E'_r = g^*(E)$
       such that $E'_{i}/E'_{i-1}$ is locally free of rank $1$ and
       $\deg((E'_{i}/E'_{i-1})_K) = 0$ for every $i$.

 \item "(c)"
       $\deg(\rest{(E'_{i}/E'_{i-1})}{C}) = 0$ for all vertical curves
       of $X'$.
\endroster
\endproclaim

\demo{Proof}
By (6) of Proposition~9.1 and
Lemma~9.6, we may assume that $E_K$ has a filtration
$0 = F_0 \subset F_1 \subset \cdots \subset F_{r-1} \subset F_r = E_K$
such that $F_i/F_{i-1}$ is locally free of rank $1$ and
$\deg(F_i/F_{i-1}) = 0$ for every $i$.
Therefore, there is a filtration of $E$ :
$0 = E_0 \subset E_1 \subset \cdots \subset E_{r-1} \subset E_r = E$
such that $E_{i}/E_{i-1}$ is torsion free and $(E_{i})_K = F_i$.
Thus by Corollary~9.4, we have our corollary.
\qed
\enddemo

Generalizing the above corollary, we would like to pose the following
question.

\definition{Question 9.8}
Let $f : X \to \Spec(O_K)$ be a regular arithmetic surface and
$E$ a semistable vector bundle on $X$. We have two questions.
\roster
\item "(1)" Is there a canonical height function on $\SStable_{X_K/K}(r, 0)$?

\item "(2)" If it exists, is there a relation
between $\achsup{2}{E}$ and the canonical height function?
\endroster
\enddefinition
\vfill\eject

\head 10. Torsion vector bundles
\endhead

In this section,
we consider conditions for the equality of $\ach_2(E, h) \leq 0$.
First of all, let us consider the case where $E$ is a line bundle.

\proclaim{Proposition 10.1}
Let $K$ be an algebraic number field and $O_K$ the ring of integers.
Let $f : X \longrightarrow \Spec(O_K)$ be a regular arithmetic surface and
$(L, h)$ a Hermitian line bundle on $X$ with $\deg(L_{\overline{\QQ}}) = 0$.
Then, $\achern{1}{L, h}^2 = 0$ if and only if
\roster
\item "i)" $h$ is Einstein-Hermitian,

\item "ii)" $\deg(\rest{L}{C}) = 0$ for all vertical curves $C$ on $X$, and

\item "iii)" $L_K$ is a torsion point of $\Pic^0(X_K)$.
\endroster
\endproclaim

\demo{Proof}
First, we assume the conditions i), ii) and iii).
By Lemma~9.2, we get
$$
 \achern{1}{L, h}^2 = -2 [K : \QQ] \operatorname{Height}([L_K]) = 0.
$$

Next we assume that $\achern{1}{L, h}^2 = 0$.
By (3) of Proposition~9.1, $h$ is Einstein-Hermitian.
Moreover, by virtue of (4) of Proposition~9.1, we have
$\deg(\rest{L}{C}) = 0$ for all vertical curves $C$ of $f$.
Hence, by Lemma~9.2, we get
$$
 \achern{1}{L, h}^2 = -2 [K : \QQ] \operatorname{Height}([L_K]),
$$
which implies $\operatorname{Height}([L_K]) = 0$.
Therefore, $L_K$ is a torsion of $\Pic^0(X_K)$.
\qed
\enddemo

\definition{Definition 10.2}
Let $M$ be a complex manifold and $E$ a rank $r$ flat vector bundle on $M$.
The vector bundle $E$ defines a representation
$$\rho : \pi_1(M) \to \GL_r(\CC)$$
of the fundamental group $\pi_1(M)$ of $M$.
$E$ is said to be {\it of torsion type} if
the image of $\rho$ is a finite group.
\enddefinition

\proclaim{Lemma 10.3}
Let $X$ be a smooth algebraic variety over $\CC$ and $E$ a flat vector
bundle on $X$. Then, $E$ is of torsion type if and only if
there is a dominant morphism of algebraic varieties $f : Y \to X$
over $\CC$ such that the composition of homomorphisms:
$$
\pi_1(Y) \to \pi_1(X) \to \GL_r(\CC)
$$
is trivial.
\endproclaim

\demo{Proof}
First we assume that $E$ is of torsion type.
Then, since $\Ker(\rho)$ has a finite index in $\pi_1(X)$,
there is a finite etale covering $f : Z \to X$ of algebraic
varieties such that $\pi_1(Z) = \Ker(\rho)$.
Thus we have the first assertion.

Next we assume that
there is a dominant morphism of algebraic varieties $f : Y \to X$
over $\CC$ such that the composition of homomorphisms
$\pi_1(Y) \to \pi_1(X) \to \GL_r(\CC)$ is trivial.
By Proposition~2.9.1 of \cite{Kol},
the image of $\pi_1(Y) \to \pi_1(X)$
has a finite index in $\pi_1(X)$.
Thus the image of $\rho$ is a finite group.
\qed
\enddemo

\proclaim{Proposition 10.4}
Let $X$ be a smooth projective curve over an algebraically closed field
$k$ with $k \subset \CC$.
Let $E$ be a rank $r$ vector bundle on $X$ such that
$E$ is semistable and $\deg(E) = 0$.
Then, the following are equivalent.
\roster
\item "(a)" There is a finite etale covering $f : Z \to X$ of $X$
over $k$ with $f^*(E) \simeq \OO_Z^{\oplus r}$.

\item "(b)" There is a surjective morphism $g : Y \to X$ of
smooth projective varieties over $k$ with $g^*(E) \simeq \OO_Y^{\oplus r}$.

\item "(c)" There is a finite etale covering $f' : Z' \to X_{\CC}$
of $X_{\CC}$ over $\CC$ such that ${f'}^*(E_{\CC}) \simeq \OO_{Z'}^{\oplus r}$.

\item "(d)" There is a surjective morphism $g' : Y' \to X_{\CC}$ of
smooth projective varieties over $\CC$
with ${g'}^*(E_{\CC}) \simeq \OO_{Y'}^{\oplus r}$.

\item "(e)" $E_{\CC}$ is flat and of torsion type.

\endroster
\endproclaim

\demo{Proof}
First of all, we prepare the following lemma.

\proclaim{Lemma 10.5}
Let $X$ be a smooth projective curve over an algebraically closed
field of characteristic zero and $E$ a rank $r$ vector bundle on $X$.
If $E$ is semistable, $\deg(E) = 0$ and $\pi^*(E) \simeq \OO_Y^{\oplus r}$
for some finite covering $\pi : Y \to X$, then $E$ is poly-stable.
\endproclaim

\demo{Proof}
We prove this lemma by induction on $\rank E$.
Clearly we may assume that $E$ is not stable. Moreover,
we may assume that $\pi : Y \to X$ is a Galois covering.
Since $E$ is not stable, there is an exact sequence:
$$
  0 \to F \to E \to Q \to 0
$$
such that $F$ and $Q$ are locally free, $\deg(F) = \deg(Q) = 0$ and
that $F$ and $Q$ are semistable.
Then, it is easy to see that the exact sequence:
$$
  0 \to \pi^*(F) \to \pi^*(E) \to \pi^*(Q) \to 0
$$
splits, $\pi^*(F) \simeq \OO_Y^{\oplus s}$ and that
$\pi^*(Q) \simeq \OO_Y^{\oplus t}$,
where $s = \rank F$ and $t = \rank Q$.
Hence $F$ and $Q$ are poly-stable by hypothesis of induction.
Therefore, it is sufficient to see that the natural homomorphism
$$
  \Ext_X^1(Q, F) \longrightarrow \Ext_Y^1(\pi^*(Q), \pi^*(E))
$$
is injective because the injectivity of the above homomorphism
implies that the exact sequence
$$
  0 \to F \to E \to Q \to 0
$$
splits.
Since $\pi$ is the Galois covering,
we have a trace map $\pi_*(\OO_Y) \to \OO_X$, which shows us that
an exact sequence:
$$
  0 \to \OO_X \to \pi_*(\OO_Y) \to \pi_*(\OO_Y)/\OO_X \to 0
$$
splits. Therefore,
$  H^1(X, F \otimes Q^{\vee}) \longrightarrow
   H^1(X, F \otimes Q^{\vee} \otimes \pi_*(\OO_Y)) $
is injective. Thus
$\Ext_X^1(Q, F) \longrightarrow \Ext_Y^1(\pi^*(Q), \pi^*(E))$
is injective.
\qed
\enddemo

\bigskip
Let us start the proof of Proposition~10.4.
(a) $\Longrightarrow$ (b),
(b) $\Longrightarrow$ (d) and
(e) $\Longrightarrow$ (c) are trivial.
So it is sufficient to show (d) $\Longrightarrow$ (e) and
(c) $\Longrightarrow$ (a).

\medskip
(d) $\Longrightarrow$ (e) :
By Lemma~10.5,  $E_{\CC}$ is poly-stable.
Thus $E_{\CC}$ is a flat vector bundle.
So $E_{\CC}$ comes from a representation
$$
 \rho : \pi_1(X_{\CC}) \to \GL_r(\CC)
$$
of the fundamental group $\pi_1(X_{\CC})$ of $X_{\CC}$.
Since ${g'}^*(E_{\CC})$ is also flat, it has an Einstein Hermitian metric
$h$. On the other hand, since ${g'}^*(E_{\CC}) \simeq \OO_{Y'}^{\oplus r}$,
by \cite{Ko, Chap. V, Proposition~8.2}, $({g'}^*(E_{\CC}), h)$ is
isometric to $(\OO_{Y'}, h_1) \oplus \cdots \oplus (\OO_{Y'}, h_r)$
as Einstein-Hermitian vector bundles, which shows us that
the composition of homomorphisms:
$$
   \pi_1(Y') \to \pi_1(X_{\CC}) \to \GL_r(\CC)
$$
is trivial. Thus, by Lemma~10.3, $E_{\CC}$ is of torsion type.

\medskip
(c) $\Longrightarrow$ (a) :
Clearly we can find a field $F$ such that
$k \subset F \subset \CC$, $F$ is finitely generated over
$k$ and that $f' : Z' \to X_{\CC}$ and
${f'}^*(E_{\CC}) \simeq \OO_{Z'}^{\oplus r}$ are defined over $F$.
Thus there are algebraic varieties $T$ and $\widetilde{Z}$ over
$k$ and projective morphisms $h : \widetilde{Z} \to T$
and $\tilde{f} : \widetilde{Z} \to X \times T$ over $k$
with the following properties.
\roster
\item "(1)" The function field of $T$ is $F$.

\item "(2)" If $p : X \times T \to T$ and $q : X \times T \to X$
are the natural projections,
then we have $p \cdot \tilde{f} = h$.

\item "(3)" $\tilde{f} \times_T \Spec(\CC) :
\widetilde{Z} \times_T \Spec(\CC) \to (X \times T) \times_T \Spec(\CC)$
is nothing more that $f' : Z' \to X_{\CC}$.
\endroster
$$
\CD
\widetilde{Z} @>{\tilde{f}}>> X \times T \\
@V{h}VV                       @VV{p}V \\
T             @=              T
\endCD
$$
At the generic point $\eta$ of $T$, we have $\tilde{f}_{\eta}$ is etale
and $(\tilde{f}^*(q^*(E)))_{\eta} \simeq
\OO_{\widetilde{Z}_{\eta}}^{\oplus r}$.
Hence there is a non-empty open set $U$ of $T$ such that
$\tilde{f}_{t}$ is etale and $(\tilde{f}^*(q^*(E)))_{t} \simeq
\OO_{\widetilde{Z}_{t}}^{\oplus r}$ for all $t \in U$.
Therefore, if we choose a closed point $t_0$ of $U$, we have (a).
\qed
\enddemo

\definition{Definition 10.6}
Let $K$ be an algebraic number field and $O_K$ the ring of integers.
Let $f : X \longrightarrow \Spec(O_K)$ be a regular arithmetic surface and
$E$ a vector bundle on $X$. $E$ is said to be {\it of torsion type}
if $E_{\overline{\QQ}}$ satisfies one of equivalent conditions
(a) -- (d) in Proposition~10.4.
\enddefinition

\definition{Question 10.7}
Let $K$ be an algebraic number field and $O_K$ the ring of integers.
Let $f : X \longrightarrow \Spec(O_K)$ be a regular arithmetic surface and
$(E, h)$ a Hermitian vector bundle on $X$ such that
$\deg(E_K) = 0$ and $E_{\overline{\QQ}}$ is semistable.
Then, if $\ach_2(E, h) = 0$, is $E$ of torsion type?
\enddefinition

For this question, we have the following partial answer.

\proclaim{Proposition 10.8}
Let $K$ be an algebraic number field and $O_K$ the ring of integers.
Let $f : X \longrightarrow \Spec(O_K)$ be a regular arithmetic surface and
$(E, h)$ a rank $r$ Hermitian vector bundle on $X$ such that
$\deg(E_K) = 0$ and $E_{\overline{\QQ}}$ is semistable.
Assume that there is a filtration of $E_{\overline{\QQ}}$:
$$
   0 = E_0 \subset E_1 \subset E_2 \subset \cdots \subset E_{r-1}
       \subset E_r = E_{\overline{\QQ}}
$$
such that $E_i/E_{i-1}$ is a locally free sheaf of rank $1$ and
$\deg(E_i/E_{i-1}) = 0$ for every $1 \leq i \leq r$.
Then, if $\ach_2(E, h) = 0$, $E$ is of torsion type.
\endproclaim

\demo{Proof}
We may assume that the filtration is defined over $K$.
Hence we can construct a filtration of $E$:
$$
   0 = F_0 \subset F_1 \subset F_2 \subset \cdots \subset F_{r-1}
       \subset F_r = E
$$
such that $\rest{F_i}{\overline{\QQ}} = E_i$,
$F_i$ is locally free and that $F_i/F_{i-1}$ is torsion free.
Let $h_{F_i}$ be the induced sub-metric of $h$ and $h_i$ the
quotient metric of $F_i/F_{i-1}$ induced by $h_{F_i}$.
Let $L_i$ be the double dual of $F_i/F_{i-1}$. Then, by Corollary~7.4,
we have
$$
\align
0 = \ach_2(E, h) &
 \leq \ach_2((F_1/F_0, h_1) \oplus \cdots \oplus (F_r/F_{r-1}, h_r)) \\
& \leq \ach_2((L_1, h_1) \oplus \cdots \oplus (L_r, h_r)) \leq 0.
\endalign
$$
Therefore, we get
$$
\ach_2((L_1, h_1) \oplus \cdots \oplus (L_r, h_r)) = \frac{1}{2}
(\achern{1}{L_1, h_1}^2 + \cdots + \achern{1}{L_r, h_r}^2) = 0
$$
and $(E, h)$ is isometric to $(L_1, h_1) \oplus \cdots \oplus (L_r, h_r)$
on each infinite fiber.
Since $\achern{1}{L_i, h_i}^2 \leq 0$ for each $i$,
we obtain $\achern{1}{L_i, h_i}^2 = 0$. Therefore, by Proposition~10.1,
$(L_i)_{\CC}$ is of torsion type. Therefore $E_{\CC}$ is torsion type
because $E_{\CC} = (L_1)_{\CC} \oplus \cdots \oplus (L_r)_{\CC}$.
\qed
\enddemo
\vfill\eject



\widestnumber\key{GS90b}
\Refs

\ref\key BC
\by R. Bott and S. S. Chern
\paper Hermitian vector bundles and the equidistribution of the zero
of their holomorphic sections
\jour Acta Math.
\vol 114
\yr 1965
\pages 71--112
\endref

\ref\key BSG 
\by J.-M. Bismut, H. Gillet and C. Soul\'{e}
\paper Analytic torsion and holomorphic determinant bundles I :
Bott-Chern forms and analytic torsion
\jour Comm. Math. Phys.
\vol 115
\yr 1988
\pages 49--78
\endref

\ref\key BV 
\by J.-M. Bismut and E. Vasserot
\paper The asymptotics of the Ray-Singer analytic torsion associated with
high power of a positive line bundle
\jour Commun. Math. Phys.
\vol 125
\yr 1989
\pages 355--367
\endref

\ref\key Bo 
\by F. A. Bogomolov
\paper Holomorphic tensors and
vector bundles on projective varieties
\jour Math. USSR-Izv.
\vol 13
\yr 1978
\pages 499--555
\endref

\ref\key Do83 
\by S. K. Donaldson
\paper A new proof of a theorem of Narasimham and Seshadri
\jour J. Diff. Geometry
\vol 18
\yr 1983
\pages 269--278
\endref

\ref\key Do85 
\by \bysame
\paper Anti self-dual Yang-Mills connections over complex algebraic surfaces
and stable bundles
\jour Proc. London Math. Soc.
\vol 50
\yr 1985
\pages 1--26
\endref

\ref\key Fa1
\by G. Faltings
\paper Calculus on arithmetic surfaces
\jour Ann. of Math.
\vol 119
\yr 1984
\pages 387--424
\endref

\ref\key Fa2 
\by G. Faltings
\book Lectures on the arithmetic Riemann-Roth Theorem
\bookinfo Annals of Mathematics Studies
\vol 127
\publ Princeton
\endref

\ref\key Gi 
\by D. Gieseker
\paper On a theorem of Bogomolov on Chern classes
of stable bundles
\jour Amer. J. Math.
\vol 101
\yr 1979
\pages 79--85
\endref

\ref\key GS88 
\by H. Gillet and C. Soul\'{e}
\paper Amplitude arithm\'{e}tique
\jour C. R. Acad. Sci. Paris
\vol t.307, S\'{e}rie I
\yr 1988
\pages 887--890
\endref

\ref\key GS90a 
\by \bysame
\paper Arithmetic Intersection Theory
\jour Publ. Math. (IHES)
\vol 72
\yr 1990
\pages 93--174
\endref

\ref\key GS90b 
\by \bysame
\paper Characteristic classes for algebraic vector bundles with
hermitian metric, I, II,
\jour Ann. of Math.
\vol 131
\yr 1990
\pages 163--203, 205--238
\endref

\ref\key GS92 
\by \bysame
\paper An arithmetic Riemann-Roch theorem
\jour Invent. math.
\vol 110
\yr 1992
\pages 473--543
\endref

\ref\key Gr
\by D. R. Grayson
\paper Reduction theory using semistability
\jour Comment. Math. Helvetici
\vol 59
\yr 1984
\pages 600--634
\endref

\ref\key Hr
\by P. Hriljac
\paper Heights and Arakelov's intersection theory
\jour Amer. J. Math.
\vol 107
\yr 1985
\pages 23--38
\endref

\ref\key Ko 
\by S. Kobayashi
\book Differential geometry of complex vector bundles
\bookinfo Publications of the Mathematical Society of Japan
\vol 15
\publ Iwanami Shoten, Publishers and Princeton University Press
\endref

\ref\key Kol
\by J. Koll\'{a}r
\paper Shafarevich maps and plurigenera of algebraic varieties
\jour preprint
\yr 1992
\endref

\ref\key Ma1
\by M. Maruyama
\paper Openness of a family of torsion free sheaves
\jour J. Math. Kyoto Univ.
\vol 16
\yr 1976
\pages 627-637
\endref

\ref\key Ma2
\by M. Maruyama
\paper Moduli of stable sheaves, I, II
\jour J. Math. Kyoto Univ.
\vol 17, 18
\yr 1977, 1978
\pages 91-126, 557-614
\endref

\ref\key MS
\by S. Mukai and F. Sakai
\paper Maximal subbundles of vector bundles on a curve
\jour manuscripta math.
\vol 52
\year 1985
\pages  251-256
\endref

\ref\key NS 
\by M. S. Narasimhan and C. S. Seshadri
\paper Stable and unitary vector bundles on compact Riemann surfaces
\jour Ann. of Math.
\vol 82
\yr 1965
\pages 540--567
\endref

\ref\key SABK
\by C. Soul\'{e}, D. Abramovich, J.-F. Burnol and J. Kramer
\book Lectures on Arakelov Geometry
\bookinfo Cambridge studies in advanced mathematics
\vol 33
\publ Cambridge University Press
\endref

\ref\key St
\by U. Stuhler
\paper Eine Bemerkung zur Reductionstheorie quadratischer Formen
\jour Arch. der Math.
\vol 27
\yr 1976
\pages 604--610
\endref

\ref\key UY 
\by K. Uhlenbeck and S.-T. Yau
\paper On the existence of Hermitian Yang-Mills connections in stable
vector bundles
\jour Comm. Pure Appl. Math.
\vol 39
\yr 1986
\pages 257--293
\endref

\ref\key Vo 
\by P. Vojta
\paper Siegel's theorem in the compact case
\jour Ann. of Math.
\vol 133
\yr 1991
\pages 509--548
\endref

\endRefs
\enddocument